\documentclass[11pt, reqno]{amsart}

\usepackage{amssymb}  

\input cyracc.def
\newfam\cyrfam
\font\tencyr=wncyr10 
\def\cyr{\fam\cyrfam\tencyr\cyracc}

\font\twelvecyr=wncyr10 scaled 1200 
\def\bigcyr{\fam\cyrfam\twelvecyr\cyracc}

\theoremstyle{plain}
\newtheorem{theorem}{Theorem}[section]
\newtheorem{lemma}[theorem]{Lemma}
\newtheorem{proposition}[theorem]{Proposition}
\newtheorem{corollary}[theorem]{Corollary}

\theoremstyle{definition}
\newtheorem{definition}[theorem]{Definition}
\theoremstyle{remark}

\newtheorem{remark}[theorem]{Remark}
\newtheorem{construction}[theorem]{Construction}

\numberwithin{equation}{section}

\newcommand{\seclabel}[1]{\label{sec:#1}} 
\newcommand{\thmlabel}[1]{\label{thm:#1}} 
\newcommand{\lemlabel}[1]{\label{lem:#1}} 
\newcommand{\corolabel}[1]{\label{coro:#1}} 
\newcommand{\proplabel}[1]{\label{prop:#1}} 
\newcommand{\deflabel}[1]{\label{def:#1}} 
\newcommand{\remlabel}[1]{\label{rem:#1}} 
\newcommand{\eqlabel}[1]{\label{eq:#1}} 
\newcommand{\conslabel}[1]{\label{cons:#1}} 

\newcommand{\secref}[1]{\ref{sec:#1}} 
\newcommand{\thmref}[1]{\ref{thm:#1}} 
\newcommand{\lemref}[1]{\ref{lem:#1}} 
\newcommand{\cororef}[1]{\ref{coro:#1}} 
\newcommand{\propref}[1]{\ref{prop:#1}} 
\newcommand{\remref}[1]{\ref{rem:#1}} 
\renewcommand{\eqref}[1]{\ref{eq:#1}} 
\newcommand{\peqref}[1]{(\eqref{#1})} 
\newcommand{\consref}[1]{\ref{cons:#1}} 

\newcommand{\Aut}{\mathrm{Aut}}
\newcommand{\End}{\mathrm{End}}

\newcommand{\MMM}{\mathcal{M}} 
\newcommand{\III}{\mathrm{I}} 

\newcommand{\sbl}[1]{\langle#1\rangle}   
\newcommand{\iv}{^{-1}} 
\newcommand{\ld}{\backslash} 
\newcommand{\rd}{/}  

\newcommand{\Id}{\mathrm{Id}}
\newcommand{\AAA}{\mathcal{A}}
\newcommand{\BBB}{\mathcal{B}}
\newcommand{\CCC}{\mathcal{C}}

\title[The Structure of F-Quasigroups]
{The Structure of F-Quasigroups}

\author[T.~Kepka]{Tom\'{a}\v{s}~Kepka$^*$}
\thanks{$^*$This work is a part of the research project MSM 0021620839
financed by MSMT and partly supported by the grant agency of the
Czech Republic, grant \# 201/05/0002}
\address{Department of Algebra \\
MFF UK, Sokolovsk\'{a} 83 \\
186 75 Praha 8, Czech Republic}
\email{kepka@karlin.mff.cuni.cz}
\author[M.~K.~Kinyon]{Michael~K.~Kinyon}
\address{Department of Mathematics \\
University of Denver \\
Denver, CO 80208 USA}
\email{mkinyon@math.du.edu}
\urladdr{http://www.math.du.edu/\symbol{126}mkinyon}
\author[J.~D.~Phillips]{J.~D.~Phillips}
\address{Department of Mathematics \& Computer Science \\
Wabash College \\
Crawfordsville, IN 47933 U.S.A.}
\email{phillipj@wabash.edu}
\urladdr{http://persweb.wabash.edu/facstaff/phillipj/}

\date{\today}


\begin{document}

\begin{abstract}
We solve a problem of Belousov which has been open since 1967: to
characterize the loop isotopes of F-quasigroups. We show that every F-quasigroup
has a Moufang loop isotope which is a central product of its nucleus and
Moufang center. We then use the loop to reveal the structure of the
associated F-quasigroup.
\end{abstract}

\maketitle


\section{Introduction}
\seclabel{intro}

A \emph{quasigroup} $(Q,\cdot)$ is a set $Q$ together with a binary
operation $\cdot : Q\times Q\to Q$ such that for each $a,b\in Q$, the
equations $ax = b$ and $ya = b$ have unique solutions $x,y\in Q$. 
Equivalently, we may consider a quasigroup
$(Q,\cdot,\ld,\rd)$ to be a set $Q$ together with three binary 
operations $\cdot,\ld,\rd : Q\times Q\to Q$ such that the equations
$x\ld (xy) = y$, $(xy)\rd y = x$, $x(x\ld y) = y$, and $(x\rd y)y = x$
hold for all $x,y\in Q$. A quasigroup with a neutral element is
called a \emph{loop}.
For the basic theory of quasigroups and loops, we refer to the standard
texts \cite{Bel,Br,Pf}.

Because we will be dealing extensively with sets which possess both a
quasigroup structure and a loop structure, we will use additive notation
for loops, even though the loop may not necessarily be commutative.
Thus the neutral element of a loop $(Q,+)$ is denoted by $0$.

Among the earliest studied varieties of quasigroups are \emph{F-quasigroups}.
These were introduced by Murdoch in 1939 \cite{Mu}.
F-quasigroups are defined by the equations
\begin{align*}
  x \cdot yz &= xy\cdot (x\ld x)z  \tag{$F_l$} \\
  zy\cdot x &= z(x\rd x) \cdot yx \tag{$F_r$} 
\end{align*}
A quasigroup satisfying ($F_l$) (resp. ($F_r$)) is said to be a
\emph{left} (resp. \emph{right}) \emph{F-quasigroup}. 
Murdoch did not actually name this particular variety. The earliest
use of the term ``F-quasigroup'' we can find is in a paper of
Belousov \cite{BO}.

Quasigroups $(Q,\cdot)$ and
$(\tilde{Q},\tilde{\cdot})$
are said to be \emph{isotopic} if there are 
bijections $f,g,h: Q\to \tilde{Q}$ such that
$f(x)\tilde{\cdot} g(y) = h(x\cdot y)$ for all $x,y\in Q$.
If $Q = \tilde{Q}$ and $h = \mathrm{id}_Q$, then $(Q,\cdot)$ and
$(\tilde{Q},\tilde{\cdot})$ are said to be \emph{principally isotopic}. Every quasigroup
$(Q,\cdot,\ld,\rd)$ is principally isotopic to a loop $(Q,+)$:
fix $a,b\in Q$, and set $0 = ba$ and $x+ y = (x\rd a)(b\ld y)$ for all
$x,y\in Q$.

One of the key issues in the study of any variety of quasigroups is 
to characterize those loops to which quasigroups in the variety are
isotopic. The exemplar of results of this type is the (weak form of the) Toyoda-Bruck
theorem: every medial quasigroup is isotopic to an abelian group
\cite{To,Br44}. This result was later generalized to 
distributive quasigroups \cite{BelDis}, trimedial quasigroups \cite{Ke76},
and semimedial quasigroups \cite{Ke78a, Ke78}.

In the very first of the Problems ({\bigcyr zadaqi}) in Belousov's 1967 book
(\cite{Bel}, p. 216), we find the following question, which has remained open until now:
\medskip

\begin{quotation}\small{
1. \ldots {\cyr Kakim lupat izotopny dvustoronnie} F-{\cyr{kvazigruppy}}\ldots ?
\medskip

\noindent 1. \ldots To which loops are two-sided F-quasigroups isotopic\ldots ?
(our translation)
}
\end{quotation}
\medskip

\noindent In 1979, one of us \cite{Ke79} implicitly conjectured the following
answer to Belousov's question.
\smallskip

\begin{center}
\emph{Every loop $(Q,+)$ isotopic to a given F-quasigroup} $(Q,\cdot)$
\emph{is a Moufang loop.}
\end{center}
\medskip

\noindent Recall that \emph{Moufang loops} $(Q,+)$ are defined by any of
the equivalent identities
\[
\begin{array}{ccc}
x + (y + (x + z)) = ((x+y)+ x)+z
&\qquad &
x+((y+z)+ x) = (x+y)+ (z+x)
\\
((z+x)+ y)+x = z+(x+ (y+x))
&\qquad &
(x+ (y+z))+ x = (x+y)+ (z+x)
\end{array}
\]
for all $x,y\in Q$ \cite{Bel,Br,Pf}. Moufang loops are isotopically invariant,
that is, every loop isotopic to a Moufang loop is a Moufang
loop. In view of this, the conjecture can be stated more concisely:
\medskip

\begin{center}
\emph{Every F-quasigroup is isotopic to a Moufang loop.}
\end{center}
\medskip

In \cite{Ke79}, it was shown that if a given F-quasigroup 
is isotopic to a Moufang loop, then at least one of the loop isotopes has
additional structure. This will be seen below in the
statement of our main result, but first, we need additional notation.

Recall that the \emph{nucleus} of a loop $(Q,+)$ is defined by
\[
\begin{array}{rcl}
N(Q) &= \{ a\in Q : & (a+x)+ y = a+ (x+y), \quad (x+a)+ y = x+ (a+y), \\ 
&& \text{and}\quad x+ (y+a) = (x+y)+ a, \quad \forall x,y\in Q\} .
\end{array}
\]
In a Moufang loop $(Q,+)$, $N(Q)$ is a normal associative subloop of
$(Q,+)$ \cite{Bel,Br,Pf}.
The \emph{Moufang center} of a loop is
\begin{eqnarray*}
K(Q) &= \{ a\in Q : (a+a)+ (x+y) = (a+x)+ (a+y), \ \ \forall x,y\in Q \} \\
&= \{ a\in Q : (x+y)+ (a+a) = (x+a)+ (y+a), \ \ \forall x,y\in Q \}.
\end{eqnarray*}
In any loop, the Moufang center is a commutative Moufang
subloop (\cite{Br}, p. 94).

\begin{definition}
\deflabel{NK-loop}
An \emph{NK-loop} is a loop $(Q,+)$ satisfying $Q = N(Q)+ K(Q)$,
that is, for each $a\in Q$, there exists $n\in N(Q)$, $k\in K(Q)$ such that $a = n+k$.
\end{definition}

In \S\secref{NK-loops}, we will study NK-loops in some detail. In
particular, every NK-loop is Moufang (Theorem \thmref{NK-Moufang}).

With these notions in hand, here is one of our main results.

\begin{theorem}
\thmlabel{main}
For a quasigroup $(Q,\cdot)$,
the following are equivalent:
\begin{enumerate}
\item[1.] $(Q,\cdot)$ is an F-quasigroup.
\item[2.] There exist a Moufang NK-loop $(Q,+)$,  
$f,g\in \Aut(Q,+)$, and $e\in N(Q,+)$ such that
\[
x\cdot y = f(x)+ e+ g(y)
\]
for all $x,y\in Q$, $fg = gf$, and
$x+ f(x), x+ g(x)\in N(Q,+)$,
$-x+ f(x),-x+ g(x)\in K(Q,+)$ for all $x\in Q$.
\end{enumerate}
\end{theorem}

Theorem \thmref{main} not only confirms the conjecture
of \cite{Ke79}, but also characterizes those Moufang loops
that can occur as loop isotopes of F-quasigroups.
In addition, the theorem generalizes the corresponding results in
the distributive \cite{BelDis} and trimedial \cite{Ke76} cases.

In the aforementioned Open Problem, Belousov also asked the same
question regarding loop isotopes for the one-sided case of left
F-quasigroups. An answer of sorts was found by Golovko, who showed
that every left F-quasigroup is isotopic to a \emph{left M-loop}
\cite{Gol}. (It would take us too far afield to include the
definition here.) This result was later included by Belousov in his
lecture notes \cite{BB}, which are not easily accessible. In a differential
geometric context, the relationship between left F-quasigroups and
left M-loops was used by Sabinin and his students to study transsymmetric
spaces. See \cite{Sab} and the complete bibliography therein for
this particular line of inquiry. This book also has the most easily
accessible proof of Golovko's result.

It turns out that the notion of M-loop is of no help in dealing with
the case of two-sided F-quasigroups. Although it follows from Golovko's
work that every F-quasigroup is isotopic to a loop which is both a
left and right M-loop, that information does not seem to be sufficient
to characterize the loop isotopes, nor does the M-loop structure seem
to help in establishing Theorem \thmref{main}.

We conclude this introduction with an outline of the sequel. In \S\secref{prelim},
we recite general preliminary results on quasigroups and loops which will be used
later in the paper. In \S\secref{Moufang}, we state necessary well-known results
about Moufang loops, and also present a few technical lemmas needed later.
In \S\secref{F-quasi}, we present basic facts about (left) F-quasigroups.
In \S\secref{NK-loops}, we examine the structure of NK-loops and show that
they are Moufang (Theorem \thmref{NK-Moufang}). We also give a sufficient
condition for a loop to be NK (Theorem \thmref{NK-char}), and this turns
out to be the main tool in showing that F-quasigroups are isotopic to NK-loops.

Our proof of Theorem \thmref{main} is really split up into the
harder implication (1)$\implies$(2) in Theorem \thmref{kepka-5.15} and the
easier converse in Proposition \propref{kepka-6.2}. The former is in 
\S\secref{linear} and the latter is in \S\secref{linear2}. In the
remainder of \S\secref{linear2} as well as section \S\secref{linear3}, we study
the structure of F-quasigroups by using their representations in terms of
NK-loops. In \S\secref{forms}, we formalize the relationship between 
(pointed) F-quasigroups and NK-loops with additional data (which we call
arithmetic forms) and show an appropriate equivalence of equational classes
(and categories) in Theorem \thmref{kepka-8.5}. Finally, in \S\secref{summary}, we
reap the rewards of this equivalence and our work in \S\secref{linear2} and
\S\secref{linear3} by presenting a summary of the structure of
F-quasigroups.

This paper redevelops some of the main results of \cite{Ke75, Ke79}, as they
are necessary to give a complete proof of Theorem \thmref{main} and its consequences.
Wherever possible, we try to give shorter, clearer proofs.

We are pleased to acknowledge the assistance of OTTER, an automated
deduction tool developed by McCune \cite{Mc}. OTTER found a proof
of Theorem \thmref{NK-char}, which we then ``humanized'' into the
form presented here. 


\section{Preliminaries}
\seclabel{prelim}

In a quasigroup $(Q,\cdot,\ld,\rd)$, it is useful to introduce notation for
local right and left neutral elements. Here we adopt the following:
\[
\alpha(x) = x\ld x \qquad \beta(x) = x\rd x ,
\]
that is $x\cdot \alpha(x) = x$ and $\beta(x)\cdot x = x$.

For $a\in Q$, left and right translations $L_a, R_a : Q\to Q$
are defined by $L_a(x) = ax$ and $R_a(x) = xa$ for $x\in Q$. 
The \emph{multiplication group} of $Q$ is the group generated by
all translations: $\MMM(Q) = \sbl{L_a, R_a : a\in Q}$.  

If $(Q,+)$ is a loop, then the stabilizer of the neutral element
$0\in Q$ is called the \emph{inner mapping group} of $Q$ and
is denoted by $\III(Q)$. A loop is called an \emph{A-loop}
if $\III(Q)\subseteq \Aut(Q,+)$. 

In \S\secref{intro}, we have already defined the nucleus and
Moufang center of a loop, which are subloops. The \emph{commutant}
or \emph{semicenter} of a loop $(Q,+)$ is the set
\[
C(Q) = \{ a\in Q : a+x = x+a\ \ \forall x\in Q\} .
\]
The commutant is not necessarily a subloop.
Finally, the \emph{center} of a loop $(Q,+)$ is defined by
\[
Z(Q) = N(Q)\cap C(Q) .
\]

\begin{lemma}
\lemlabel{K_in_C}
In a loop $(Q,+)$,
\begin{enumerate}
\item $K(Q)\subseteq C(Q)$.
\item $C(Q)$ is a characteristic subset.
\item $K(Q)$ is a characteristic subloop.
\item $N(Q)$ is a characteristic subgroup.
\item $Z(Q) = N(Q)\cap K(Q)$.
\item $Z(Q)$ is a normal abelian subgroup.
\end{enumerate}
\end{lemma}

\begin{proof}
For (1): If $a\in K(Q)$, then for all $x\in Q$,
$a + (a + x) = (a + 0) + (a + x) = 2a + x = (a + x) + (a + 0) = (a + x) + a$.
Replacing $a + x$ with $x$, we have the desired result. The rest are clear
or can be found in the standard references \cite{Bel, Br, Pf}.
\end{proof}

\begin{lemma}
\lemlabel{A-loops}
Let $(Q,+)$ be an A-loop.
\begin{enumerate}
\item Every characteristic subloop is normal.
\item $N(Q)$ is a normal subloop.
\item $K(Q)$ is a normal subloop.
\end{enumerate}
\end{lemma}

\begin{proof}
Recall that a subloop is normal if and only if it is invariant under the
action of $\III(Q)$. Then (1) follows immediately, and (2) and (3) follow
from \lemref{K_in_C}.
\end{proof}

For a quasigroup $(Q,\cdot)$, we define
\[
M(Q) = \{ a\in Q : xa\cdot yx = xy\cdot ax \ \forall x,y\in Q \}
\]

\begin{lemma}
\lemlabel{M_in_C}
For a loop $(Q,+)$,
\begin{enumerate}
\item $M(Q)\subseteq C(Q)$.
\item $Z(Q) = N(Q)\cap M(Q)$.
\end{enumerate}
\end{lemma}

\begin{proof}
For $a\in M(Q)$, $(x + a) + (y + x) = (x + y) + (a + x)$ for all $x,y\in Q$.
Take $x = 0$ to obtain (1). Then (2) is clear.
\end{proof}

A quasigroup $(Q,\cdot)$ is called 
\begin{itemize}
\item \emph{medial} if $xa \cdot by = xb\cdot ay$ for all $x,y,a,b\in Q$;
\item \emph{monomedial} (\emph{dimedial}, \emph{trimedial}, resp.) if every
(at most) one-generated (two-generated, three-generated, resp.) subquasigroup
of $Q$ is medial;
\item \emph{distributive} if $x\cdot yz = xy\cdot xz$ and $zy\cdot x = zx\cdot yx$
for all $x,y,z\in Q$;
\item \emph{symmetric} if $xy = yx$ and $x\cdot xy = y$ for all $x,y\in Q$.
\end{itemize}

\begin{proposition}[\cite{To,Br44}]
\proplabel{kepka-2.7}
For a quasigroup $(Q,\cdot)$,
the following are equivalent:
\begin{enumerate}
\item[1.] $(Q,\cdot)$ is medial.
\item[2.] There exist an abelian group $(Q,+)$,  
$f,g\in \Aut(Q,+)$, and $e\in Q$ such that
$x\cdot y = f(x)+ e + g(y)$
for all $x,y\in Q$, and $fg = gf$.
\end{enumerate}
\end{proposition}

\begin{proposition}[\cite{Ke76}]
\proplabel{kepka-2.8}
For a quasigroup $(Q,\cdot)$,
the following are equivalent:
\begin{enumerate}
\item[1.] $(Q,\cdot)$ is trimedial.
\item[2.] There exist a commutative Moufang loop $(Q,+)$,  
$f,g\in \Aut(Q,+)$, and $e\in Z(Q,+)$ such that
$x\cdot y = f(x)+ e+ g(y)$
for all $x,y\in Q$, $fg = gf$, and $x + f(x),x + g(x)\in
Z(Q,+)$ for all $x\in Q$.
\end{enumerate}
\end{proposition}

\begin{proposition}[\cite{BelDis}]
\proplabel{kepka-2.11}
For a quasigroup $(Q,\cdot)$,
the following are equivalent:
\begin{enumerate}
\item[1.] $(Q,\cdot)$ is symmetric and distributive.
\item[2.] There exist a commutative Moufang loop $(Q,+)$
of exponent $3$ such that $x\cdot y = - x - y$
for all $x,y\in Q$.
\end{enumerate}
\end{proposition}

\begin{proposition}[\cite{BelDis}]
\proplabel{kepka-2.10a}
For a quasigroup $(Q,\cdot)$,
the following are equivalent:
\begin{enumerate}
\item[1.] $(Q,\cdot)$ is distributive.
\item[2.] There exist a commutative Moufang loop $(Q,+)$,
$f\in \Aut(Q,+)$ such that $x\cdot y = f(x) + y - f(y)$
for all $x,y\in Q$, $x\mapsto x - f(x)$ is a permutation
of $Q$, and $x + f(x)\in Z(Q,+)$ for all $x\in Q$.
\end{enumerate}
\end{proposition}

We conclude this section by introducing some useful groups of
pairs of mappings. For a quasigroup $(Q,\cdot)$, define
\begin{align*}
\AAA(Q) &= \{ (p,q) : Q\times Q\to Q\times Q\ |\ p(xy) = q(x)y\ \forall x,y\in Q \} \\
\BBB(Q) &= \{ (p,q) : Q\times Q\to Q\times Q\ |\ p(xy) = xq(y)\ \forall x,y\in Q \} \\
\CCC(Q) &= \{ (p,q) : Q\times Q\to Q\times Q\ |\ p(x)y = xq(y)\ \forall x,y\in Q \} \\
\end{align*}
Let
\[
\AAA_l(Q) = \{ p : (p,q)\in \AAA(Q) \}  \qquad \text{and}\qquad 
\AAA_r(Q) = \{ q : (p,q)\in \AAA(Q) \}
\]
and similarly define $\BBB_l(Q)$, $\BBB_r(Q)$, $\CCC_l(Q)$, and $\CCC_r(Q)$.
It is easy to see that every mapping from 
$\AAA_l(Q) \cup \AAA_r(Q)\cup \BBB_l(Q) \cup \BBB_r(Q) \cup \CCC_l(Q) \cup \CCC_r(Q)$
is a permutation of $Q$. These mappings are known as \emph{regular permutations}
of the quasigroup $Q$. In addition, $\AAA_l(Q) \cong \AAA_r(Q)$, $\BBB_l(Q)\cong \BBB_r(Q)$, 
$\CCC_l(Q)\cong \CCC_r(Q)$ are permutation groups.


\section{F-quasigroups}
\seclabel{F-quasi}

For convenience, we repeat here the basic definitions of \S\secref{intro},
using the notational conventions of \S\secref{prelim}. 
A quasigroup $(Q,\cdot)$ is said to be a \emph{left F-quasigroup}
if it satisfies the identity
\[
x\cdot yz = xy\cdot \alpha(x)z \tag{$F_l$}
\]
for all $x,y,z\in Q$. $(Q,\cdot)$ is said to be a \emph{right F-quasigroup} if it
satisfies
\[
zy\cdot x = z\beta(x)\cdot yx \tag{$F_r$}
\]
for all $x,y,z\in Q$. If $(Q,\cdot)$ is both a left F-quasigroup and a right F-quasigroup,
then $(Q,\cdot)$ is called a (two-sided) F-quasigroup.

\begin{lemma}
\lemlabel{kepka-2.2}
The following conditions are equivalent for a quasigroup $(Q,\cdot)$.
\begin{enumerate}
\item $(Q,\cdot)$ is a left F-quasigroup.
\item $L_x y\cdot L_{\alpha(x)}z = L_x (yz)$, $\forall x,y,z\in Q$.
\item $L_x L_y = L_{xy} L_{\alpha(x)}$, $\forall x,y\in Q$.
\item $L(x,y) = L_{\alpha(x)}$, $\forall x,y\in Q$.
\item $L(x,y) = L(x,z)$, $\forall x,y,z\in Q$.
\item $L_x R_z = R_{\alpha(x)z} L_x$, $\forall x,z\in Q$.
\end{enumerate}
\end{lemma}

\begin{proof}
These are all just different ways of rewriting the definition of left F-quasigroup.
\end{proof}

\begin{lemma}
\lemlabel{kepka-2.3}
Let $(Q,\cdot)$ be a left F-quasigroup. Then
\begin{enumerate}
\item $\alpha \beta = \beta \alpha$ and $\alpha\in \End(Q,\cdot)$.
\item $R_a L_b = L_b R_a$ for $a,b\in Q$ if and only if $\alpha(b) = \beta(a)$.
\item $R_{\alpha(a)} L_{\beta(a)} = L_{\beta(a)} R_{\alpha(a)}$.
\end{enumerate}
\end{lemma}

\begin{proof}
(i) We compute
\[
x\cdot \alpha\beta(x)\alpha(x) = \beta(x)x \cdot \alpha\beta(x)\alpha(x)
= \beta(x)\cdot x\alpha(x) = \beta(x)x = x = x\alpha(x) .
\]
Cancelling $x$ and then dividing on the right by $\alpha(x)$, we obtain
$\alpha\beta(x) = \beta\alpha(x)$, as claimed. Further,
$xy\cdot \alpha(x)\alpha(y) = x\cdot y\alpha(y) = xy = xy\cdot \alpha(xy)$
and so $\alpha(x)\alpha(y) = \alpha(xy)$ for all $x,y\in Q$.

(ii) Observe that
\[
b\cdot xa = bx \cdot \alpha(b)a \qquad \text{and} \qquad
bx\cdot a = bx \cdot \beta(a)a ,
\]
for all $x,a,b\in Q$. The desired result follows immediately.

Finally, (iii) follows from (i) and (ii).
\end{proof}

\begin{corollary}
\corolabel{kepka-2.4}
If $(Q,\cdot)$ is an F-quasigroup, then $\alpha$ and $\beta$ are commuting
endomorphisms of $Q$.
\end{corollary}

\begin{lemma}
\lemlabel{kepka-2.5}
A loop $(Q,+)$ is a left (right, two-sided) F-loop if and only if it is a group.
\end{lemma}

\begin{proof}
This is immediate from $(F_l)$.
\end{proof}

\begin{remark}
\remlabel{kepka-2.9}
By Lemma \lemref{kepka-2.5}, every group is an F-quasigroup. In addition,
it is clear from definitions that every trimedial quasigroup is an F-quasigroup.
Thus we observe that the variety of F-quasigroups is strictly larger than
the variety of trimedial quasigroups.
\end{remark}

Again just by comparing definitions, we have the following.

\begin{proposition}
\proplabel{kepka-2.10b}
For a quasigroup $(Q,\cdot)$,
the following are equivalent:
\begin{enumerate}
\item[1.] $(Q,\cdot)$ is distributive.
\item[2.] $(Q,\cdot)$ is idempotent and trimedial.
\item[3.] $(Q,\cdot)$ is an idempotent F-quasigroup.
\end{enumerate}
\end{proposition}


\section{Moufang loops}
\seclabel{Moufang}

We begin by reciting some basic results, most of which are in the
literature.

\begin{proposition}
\proplabel{Mfg_basic}
Let $(Q,+)$ be a Moufang loop. Then
\begin{enumerate}
\item[1.] $(Q,+)$ is \emph{diassociative}, that is, for each $x,y\in Q$, $\sbl{x,y}$ is a group.
\item[2.] $N(Q)$ is a normal subloop.
\item[3.] $K(Q) = C(Q) = M(Q)$ is a characteristic subloop.
\item[4.] If $K(Q)$ is a normal subloop, then for each $x\in Q$,
the subloop generated by $\{x\}\cup K(Q)$ is commutative.
\end{enumerate}
\end{proposition}

\begin{proof}
Most of these are standard facts; see, for instance, \cite{Bel,Br,Pf}.
The only perhaps unfamiliar assertions are the second equality of (3)
and (4).

For (3): For $a\in C(Q)$, $x,y\in Q$, note that $(x + a) + (y + x) =
(x + (a + y)) + x = (x + (y + a)) + x = (x + y) + (a + x)$. Thus $c\in M(Q)$.
The other inclusion follows from Lemma \lemref{M_in_C}(1).

For (4): Let $P = \sbl{x,K(Q)}$. Then $K(Q)\subseteq K(P)$, and
it suffices to show that $x\in K(P)$. Fix $y\in P$ and set $R = \sbl{x,y}$.
Then $R$ is a group and $R / (R\cap K(Q))$ is isomorphic to a subgroup
of the cyclic group $P/K(Q)$. Thus $R/Z(R)$ is cyclic, so $R = Z(R)$,
that is, $R$ is abelian and $x+y = y+x$.
\end{proof}

Recall that a permutation $f :Q\to Q$ of a loop $(Q,+)$ is a 
left \emph{pseudoautomorphism} with left \emph{companion} $c\in Q$
if $c+ f(x+ y) = (c + f(x)) + f(y)$ for all $x,y\in Q$. Right
pseudoautomorphisms with their right companions are defined similarly.

\begin{lemma}
\lemlabel{pseudo}
Let $(Q,+)$ be a Moufang loop.
\begin{enumerate}
\item A permutation $f$ of $Q$ is a left pseudoautomorphism with left companion $c\in Q$,
iff $f$ is a right pseudoautomorphism with right companion $-c\in Q$.
\item Every inner mapping is a (left and right) pseudoautomorphism.
\item If $f$ is a pseudoautomorphism, then for each $a\in N(Q)$
and all $x\in Q$, $f(a+x) = f(a)+f(x)$, $f(x+a) = f(x)+f(a)$, and 
$f(a)\in N(Q)$.
\item A pseudoautomorphism $f$ with companion $c$ is an automorphism if and only
if $c\in N(Q)$.
\end{enumerate}
\end{lemma}

\begin{proof}
We prove (1) and leave (2) and (3) to the references, noting that (4) is clear.
If $f$ is a left pseudautomorphism with left companion $c$, then
\begin{align*}
f(x+y) - c &= -[c - f(x+y)] = -[c + f(-y-x)] = -[(c + f(-y)) + f(-x)]\\
&= f(x) - (c + f(-y)) = f(x) + (f(y) - c)
\end{align*}
for $x,y\in Q$.
\end{proof}

\begin{lemma}
\lemlabel{A-loop}
Let $(Q,+)$ be a Moufang A-loop. Then
\begin{enumerate}
\item[1.] $K(Q)$ is a normal subloop.
\item[2.] For each $x\in Q$, $3x \in N(Q)$.
\end{enumerate}
\end{lemma}

\begin{proof}
(1) is from Lemma \lemref{A-loops}. For (2),
just note that the inner mapping $y \mapsto x + (y - x)$ is both a 
pseudoautomorphism with companion $3x$ (since $Q$ is Moufang)
and an automorphism (since $Q$ is an A-loop), and apply 
Lemma \lemref{pseudo}(4).
\end{proof}

\begin{lemma}
\lemlabel{aut_C}
In a Moufang loop $(Q,+)$, the following conditions are equivalent for 
an automorphism $f$ of $Q$:
\begin{enumerate}
\item[1.] $-x + f(x)\in C(Q)$, $\forall x\in Q$
\item[2.] $f(x) - x \in C(Q)$, $\forall x\in Q$
\item[3.] $-x+ f(x) = f(x)- x \in C(Q)$, $\forall x\in Q$
\item[4.] $-x + f\iv(x)\in C(Q)$, $\forall x\in Q$
\item[5.] $f\iv(x)- x \in C(Q)$, $\forall x\in Q$
\item[6.] $-x+ f\iv(x) = f\iv(x) - x \in C(Q)$, $\forall x\in Q$
\end{enumerate}
\end{lemma}

\begin{proof}
Trivially, (3) implies (1) and (2). If (1) holds, then $f(x) = x+ (-x+ f(x))
= (-x+ f(x))+ x$, so that $f(x)-x = x- f(x)$, \emph{i.e.}, (3) holds.
Similarly, (2) implies (3). The equivalence of (4), (5), and (6) follows from
replacing $f$ with $f\iv$. Finally, apply $f\iv$ to (3) and replace
$x$ with $-x$ to get (6), using that $C(Q)$ is a characteristic subloop.
Similarly, (6) implies (3).
\end{proof}

\begin{lemma}
\lemlabel{aut_N}
In a Moufang loop $(Q,+)$, the following conditions are equivalent for 
an automorphism $f$ of $Q$ and for a fixed $i\in \{0,1\}$:

\begin{tabular}{rlcrl}
\emph{1.} & $(-1)^i x+ f(x)\in N(Q)$, $\forall x\in Q$   & & 
\emph{3.} & $(-1)^i x+ f\iv(x)\in N(Q)$, $\forall x\in Q$ \\
\emph{2.} & $f(x) + (-1)^i x \in N(Q)$, $\forall x\in Q$ & &
\emph{4.} & $f\iv(x) + (-1)^i x \in N(Q)$, $\forall x\in Q$
\end{tabular}
\end{lemma}

\begin{proof}
If $a+b\in N(Q)$, which is a normal subloop, then $b+a = -a+ (a+b)+ a \in N(Q)$.
Thus (1) is equivalent to (2), and (3) is equivalent to (4).
Apply $f\iv$ to (1), and if $i=1$, replace $x$ with $-x$ to get (4), and
conversely, (4) implies (1).
\end{proof}

\begin{lemma}
\lemlabel{aut_N_2}
In a Moufang loop $(Q,+)$, the following conditions are equivalent for 
an automorphism $f$ of $Q$:

\begin{tabular}{rlcrl}
\emph{1.} & $-x+ f(2x)\in N(Q)$, $\forall x\in Q$ & &
\emph{3.} & $-2x + f\iv(x)\in N(Q)$, $\forall x\in Q$ \\
\emph{2.} & $f(2x) - x \in N(Q)$, $\forall x\in Q$ & &
\emph{4.} & $f\iv(x) - 2x \in N(Q)$, $\forall x\in Q$
\end{tabular}

\noindent Moreover, if $(Q,+)$ is also an A-loop, then these conditions are equivalent to
the conditions of Lemma \lemref{aut_N} with $i = 0$.
\end{lemma}

\begin{proof}
The equivalence of the conditions is proven similarly to Lemma \lemref{aut_N}. Now if $(Q,+)$
is an A-loop, then $3x\in N(Q)$ for all $x\in Q$. Thus if the conditions hold, 
$x + f(x) = x + f(-2x)+ f(3x) \in N(Q)$, using (1) and the fact that $N(Q)$ is a characteristic
subloop. But then condition (1) of Lemma \lemref{aut_N} holds. The converse is similar.
\end{proof}

\section{NK-loops}
\seclabel{NK-loops}

\begin{lemma}[Pflugfelder's Theorem]
\lemlabel{pflug}
Let $(Q,+)$ be a loop and $A : Q\to Q$ a mapping. The following
are equivalent.
\begin{itemize}
\item[1.] \qquad $(x+y)+ (z+A(x)) = x+((y+z)+ A(x))$ for all $x,y,z\in Q$,
\item[2.] \qquad $(x+y)+ (z+A(x)) = (x+ (y+z))+ A(x)$ for all $x,y,z\in Q$,
\item[3.] \qquad $(Q,+)$ is a Moufang loop and $-x+A(x)\in N(Q)$
for all $x\in Q$.
\end{itemize}
Moreover, if any (and hence all) of these conditions hold, then the subloop
$\sbl{x,A(x),y}$ is a group for each $x,y\in Q$. In addition,
\[
K(Q) = \{ a\in Q : (x + a) + (y + A(x)) = (x + y) + (a + A(x)) \ \forall x,y\in Q\}.
\]
\end{lemma}

\begin{proof}
Everything except the final assertion can be found in \cite{Pf70}.
Now if $a\in K(Q)$, then 
$(x + a) + (y + A(x)) = x + ((a + y) + A(x)) = x + ((y + a) + A(x)) = (x + y) + (a + A(x))$.
Conversely, if $a\in Q$ is such that $(x + a) + (y + A(x)) = (x + y) + (a + A(x))$ holds
for all $x,y\in Q$, then taking $x = 0$, we have $a + (y + A(0)) = y + (a + A(0))$.
But $A(0) = -0 + A(0) \in N(Q)$, and so we may cancel to obtain $a + y = y + a$, that
is, $a\in K(Q)$.
\end{proof}

\begin{theorem}
\thmlabel{NK-char}
Let $(Q,+)$ be a loop and $A : Q\to Q$ a mapping. The following
are equivalent.
\begin{itemize}
\item[1.] \quad $(x+A(x))+ (y+z) = (x+y)+ (A(x)+z)$ for all $x,y,z\in Q$,
\item[2.] \quad $(Q,+)$ is a Moufang loop, and $A(x)\in K(Q)$,
$-x+ A(x)\in N(Q)$ for all $x\in Q$.
\end{itemize}
Moreover, if either (and hence both) of these conditions hold, $(Q,+)$
is an NK-loop.
\end{theorem}

\begin{proof}
(1)$\implies$(2): We let $-x$ denote the right inverse of $x$ in $Q$,
that is, $x + (-x) = 0$. Taking $y = 0$ in (1), we have $(x + A(x)) + y = 
x + (A(x) + y)$. Thus $x + (A(x) + (y + z)) = (x + y) + (A(x) + z)$.
Taking $z = -A(x)$ and cancelling, we obtain 
$A(x) + (y + (-A(x))) = y$. Replace $y$ with $A(x) + y$ and cancel
to get 
\begin{equation}\eqlabel{NK-char}
(A(x) + y) + (-A(x)) = y .
\end{equation}
Next, $[(x + y) + A(x +  y)] + [A(x) + (-A(x + y))] =$
\begin{align*}
(x + y) + A(x) 
&= (x + y) + (A(x) + 0)
= (x + A(x)) + y \\
&= (x + A(x)) + [(A(x + y) + y) + (-A(x + y))] \\
&= [x + (A(x + y) + y)] + [A(x) + (-A(x + y))]
\end{align*}
using (1) twice, then \peqref{NK-char}, and then (1) again.
Cancelling, we obtain
\begin{align*}
(x + y) + A(x +  y) &= x + (A(x + y) + y) = [(x + y) + u] + (A(x + y) + y) \\
&= [(x + y) + A(x + y)] + (u + y),
\end{align*}
where $u = u(x,y)$ satisfies $(x + y) + u = x$. Cancelling, we have
$0 = u + y$, whence $u$ is independent of $x$ and $y = -u$. Thus
$(x + (-u)) + u = x$ for all $x,u\in Q$. Now adding $A(x)$ on the 
right of \peqref{NK-char}, we obtain $A(x) + y = y + A(x)$, i.e.,
$A(x)\in C(Q)$ for all $x\in Q$. Finally,
\begin{align*}
x + ((y + z)+ A(x)) &= x+(A(x)+ (y+z)) = (x+A(x))+ (y+z) \\
&= (x+y)+ (A(x)+z) = (x+y)+ (z+A(x)).
\end{align*}
Applying Lemma \lemref{pflug}, we have that $(Q,+)$ is Moufang and $-x+ A(x)\in N(Q)$.
Since $K(Q)=C(Q)$ in Moufang loops, the proof is complete.

(2)$\implies$(1): Using Lemma \lemref{pflug} and $A(x)\in K(Q) = C(Q)$,
$(x+y)+ (A(x)+z) = (x+y)+ (z+A(x)) = x+((y+z)+ A(x)) = x+(A(x)+ (y+z))$. Since
$\sbl{x,A(x),yz}$ is a group, $(x+y)+ (A(x)+z) = (x+A(x))+ (y+z)$.

Proof of ``Moreover'': $x = A(x) - (-x + A(x))$.

\end{proof}

\begin{lemma}
\lemlabel{nk-hom}
Let $(Q,+)$ be an NK-loop. Then the mapping $N(Q)\times K(Q)\to Q ;
(n,k)\mapsto n+k$ is an epimorphism of loops.
\end{lemma}

\begin{proof}
For $n_1,n_2\in N(Q)$, $k_1,k_2\in K(Q)$, $(n_1 + k_1)+ (n_2 +k_2) =
n_1 +(k_1 + (n_2 +k_2)) = n_1 +((n_2 + k_2) + k_1) = (n_1 +n_2)+ (k_2 +k_1)
= (n_1 +n_2)+ (k_1 +k_2)$.
\end{proof}

\begin{theorem}
\thmlabel{NK-Moufang}
Every NK-loop is a Moufang A-loop.
\end{theorem}

\begin{proof}
The class of Moufang A-loops is a variety (i.e., equational class), and hence
is closed under direct products and homomorphic images. 
Obviously groups are Moufang A-loops, and commutative Moufang loops are also
A-loops (\cite{Br}, Lemma VII.3.3, p. 116). Thus if $Q$ is an NK-loop, then
$N(Q)\times K(Q)$ is a Moufang A-loop, and by Lemma \lemref{nk-hom}, so is $Q$.
\end{proof}

\begin{corollary}
\corolabel{NK-twoway}
A loop $(Q,+)$ is an NK-loop if and only if there exists a mapping
$A: Q\to Q$ such that $(x + A(x)) + (y + z) = (x + y) + (A(x) + z)$
for all $x,y,z\in Q$.
\end{corollary}

\begin{proof}
If $(Q,+)$ is an NK-loop, then $(Q,+)$ is a Moufang loop (Theorem \thmref{NK-Moufang}),
and there exists a mapping $A: Q\to Q$ such that $A(x)\in K(Q)$ and 
$-x + A(x)\in N(Q)$ for all $x\in Q$. We may then apply Theorem \thmref{NK-char}.
The converse follows directly from that same theorem.
\end{proof}

\begin{corollary}
\corolabel{NK-struct}
Let $Q$ be an NK-loop.
\begin{enumerate}
\item $Z(Q) = Z(N(Q)) = Z(K(Q))$.
\item $N(Q)$ is a normal subgroup of $Q$, $Q/N(Q)$ is a commutative Moufang loop of exponent 3,
and $Q/N(Q) \cong K(Q)/Z(Q)$.
\item $K(Q)$ is a normal commutative subloop of $Q$, $Q/K(Q)$ is a group, and $Q/K(Q) \cong N(Q)/Z(Q)$.
\end{enumerate}
\end{corollary}

\begin{remark}
Suppose $(Q,\cdot)$ is an NK-loop and set $P = N(Q)\times K(Q)$. As noted in the proof of
Theorem \thmref{NK-Moufang}, $P$ is an NK-loop. In addition, $N(P) = N(Q)\times Z(K(Q))$,
$K(P) = Z(N(Q))\times K(Q)$, and $Z(P) = Z(N(Q))\times Z(K(Q))$. The epimorphism
$\pi : P\to Q ; (n,k)\mapsto nk$ (Lemma \lemref{nk-hom}) has kernel
$\ker \pi = \{ (x,x\iv): x\in Z(Q)\} \subseteq Z(P)$. Now $P\cong Q$ iff 
$\ker \pi = \{ (1,1)\}$ iff $Z(Q) = \{1\}$. But then $Z(P) = \{1\}$, and conversely,
if $Z(P) = \{1\}$, $P\cong Q$. In this case, $N(Q)$ is a group with trivial center
and $K(Q)$ is a commutative Moufang loop with trivial center. Note that by the Bruck-Slaby
theorem, it follows that $K(Q)$ must be infinitely generated (\cite{Br}, Chap. VIII).
\end{remark}

\begin{lemma}
\lemlabel{3.12-prelim}
Let $(Q,+)$ be a commutative Moufang loop. If $f$ is a pseudoautomorphism with
companion $c$, then $f$ is an automorphism and $c\in Z(Q)$.
\end{lemma}

\begin{proof}
We show this for $f$ a left pseudoautomorphism, the right case being dual.
For all $x,y\in Q$, $(c+f(x))+f(y) = c+f(x+y) = c + f(y+x) = (c + f(y)) + f(x)$.
But then $c\in N(Q) = Z(Q)$.
\end{proof}

\begin{lemma}
\lemlabel{kepka-3.12}
Let $(Q,+)$ be an NK-loop. If $f$ is a pseudoautomorphism with companion $c$,
then $f$ is an automorphism and $c\in N(Q)$.
\end{lemma}

\begin{proof}
First write $c = a + b$ where $a\in N(Q)$, $b\in K(Q)$, and we see easily that
$b$ is a companion of $f$.  Now for $u\in K(Q)$, $v\in N(Q)$, we use 
Proposition \propref{Mfg_basic}(6) to compute 
$f(u) + v = f(u) + f(f\iv(v)) = f(u + f\iv(v)) = f(f\iv(v) + u) = v + f(u)$.
Since $(Q,+)$ is an NK-loop, $f(u)\in K(Q)$. Similarly, $f\iv(u)\in K(Q)$,
and so (the restriction of) $f$ is a pseudoautomorphism of $K(Q)$ with
companion $b$. By Lemma \lemref{3.12-prelim} and Corollary \cororef{NK-struct}(1),
$b\in Z(K(Q)) = Z(Q)$, and so $f$ is an automorphism of $(Q,+)$.
\end{proof}

\begin{lemma}
\lemlabel{kepka-3.13}
Let $(Q,+)$ be an NK-loop. Then:
\begin{enumerate}
\item For each $x\in Q$, the subloop $\sbl{x,K(Q)}$ is commutative.
\item For each $x,y\in Q$, the subloop $\sbl{x,y,N(Q)}$ is a group.
\end{enumerate}
\end{lemma}

\begin{proof}
For (1): combine Corollary \cororef{NK-struct}(3) and Proposition
\propref{Mfg_basic}(4).

For (2): By Lemma \lemref{nk-hom}, we may assume without loss that
$Q$ is commutative. But then $N(Q) = Z(Q)$, and the assertion 
follows from the diassociativity of $Q$.
\end{proof}

\begin{lemma}
\lemlabel{kepka-3.14}
Let $(Q,+)$ be a Moufang loop, and let $f\in \Aut(Q)$ satisfy
$-x + f(x) \in K(Q)$ and $-x + f(2x)\in N(Q)$ for all $x\in Q$.
Then $(Q,+)$ is an NK-loop and
\[
x + (y + z) = (f(x) + y) + ((f(-x) + x) + z)
\]
\end{lemma}

\begin{proof}
Set $k(x) = -x + f\iv(x)$. 
By Lemmas \lemref{aut_C} and \lemref{aut_N_2}, $k(x)\in K(Q)$
and $-k(x) + x = -f\iv(x) + 2x \in N(Q)$. Since
$x = k(x) + (-k(x) + x) \in K(Q) + N(Q)$, $Q$ is an NK-loop.
By Theorem \thmref{NK-char}, $(x + k(x)) + (y + z) =
(x + y) + (k(x) + z)$. That is,
$f\iv(x) + (y + z) = (x + y) + ((-x + f\iv(x)) + z)$, and so
replacing $x$ with $f(x)$, we have the rest of the result.
\end{proof}

\begin{lemma}
\lemlabel{kepka-3.15}
Let $(Q,+)$ be a Moufang loop, and let $f,g\in \Aut(Q)$
satisfy $-x + f(2x),$ $-x + g(2x)\in N(Q)$ for all $x\in Q$.
Then
\[
K(Q) = \{ a\in Q : (f^2(x) + a) + (y + g^2(x)) = 
(f^2(x) + y) + (a + g^2(x)) \ \forall x,y\in Q \}.
\]
\end{lemma}

\begin{proof}
By Lemma \lemref{aut_N_2}, $-2x + f\iv(x)\in N(Q)$.
Since $N(Q)$ is characteristic, $g(-2x) + gf\iv(x)\in N(Q)$. 
Now $g(-2x) + gf\iv(x) = ((g(-2x) + x) - x) + gf\iv(x)$ and
since $g(-2x) + x\in N(Q)$ (again by Lemma \lemref{aut_N_2}),
we have $-x + gf\iv(x)\in N(Q)$. Thus $g(-x) + g^2f\iv(x)\in N(Q)$
since $N(Q)$ is characteristic. Setting $x = f\iv(u)$, we have
$gf\iv(-u) + g^2f^{-2}(u) \in N(Q)$. Now 
$gf\iv(-u) + g^2f^{-2}(u) = ((-u + gf\iv(u))\iv - u) + g^2f^{-2}(u)$,
and since $-u + gf\iv(u)\in N(Q)$ as before, we have 
$-u + g^2f^{-2}(u)\in N(Q)$ for all $u\in Q$. Now Pflugfelder's
Theorem applies with $A = g^2f^{-2}$ and gives that $a\in K(Q)$
if and only if $(x + a) + (y + g^2f^{-2}(x)) = 
(x + y) + (a + g^2f^{-2}(x))$ for all $x,y\in Q$. Replacing 
$x$ with $f^2(x)$, we have the desired result.
\end{proof}


\section{F-quasigroups are linear over NK-loops}
\seclabel{linear}

Throughout this section, let $(Q,\cdot)$ be an F-quasigroup and let
$a,b\in Q$ be such that $hk = kh$ where $h = R_a$ and $k = L_b$
(see Lemma \lemref{kepka-2.3}). Further, put
\[
f = h R_{\beta(a)}h\iv, \quad g = k L_{\alpha(b)} k\iv, \quad
p = hk\alpha h\iv, \quad q = kh\beta k\iv (= hk\beta k\iv ) .
\]
Observe that the mappings $f,g,p,q : Q\to Q$ are permutations of $Q$.
Finally, put
\[
x + y = h\iv(x) \cdot k\iv(y)
\]
for all $x,y\in Q$ and set $0 = ba$. Then $(Q,+)$ is a (possibly
noncommutative) loop isotopic to $(Q,\cdot)$ and $0$ is the neutral element.
Our goals are to show that $(Q,+)$ is an NK-loop and that $(Q,\cdot)$
is linear over the loop. To get there, we need a sequence of lemmas.

\begin{lemma}
\lemlabel{kepka-5.2-A}
\begin{enumerate}
\item $h(x + y) = f(x) + h(y)$ for all $x,y\in Q$.
\item $h(x) = f(x) + c$ for all $x\in Q$, where $c = h(0) = ba\cdot a$.
\item $f$ is a right pseudoautomorphism with right companion $c$ of $(Q,+)$.
\end{enumerate}
\end{lemma}

\begin{proof}
For (1): We have
\[
h(h(u) + k(v)) + k(w) = uv\cdot w = u\beta(w)\cdot vw 
= h(h(u) + k\beta(w)) + k(h(v) + k(w))
\]
for all $u,v,w\in Q$ by $(F_r)$. Setting $x = h(u)$,
$y = k(v)$, and $z = k(w)$, we get
$h(x + y) + z = h(x + k\beta k\iv(z)) + k(hk\iv(y) + z)$
for all $x,y,z\in Q$. In particular, for $z = 0$, we get
$h(x + y) = h(x + k\beta k\iv(0)) + h(y)$, since $khk\iv = h$. 
Moreover, $k\beta k\iv(0) = k\beta k\iv(ba) = k\beta(a)$
and $h(x + k\beta(a)) = h(h\iv(x)\cdot \beta(a)) = f(x)$.
Thus $h(x + y) = f(x) + h(y)$ for all $x,y\in Q$ as claimed.

(2) follows from (1) by taking $y=0$, and (3) follows from (1)
and (2).
\end{proof}

\begin{lemma}
\lemlabel{kepka-5.2-B}
\begin{enumerate}
\item $k(x + y) = k(x) + g(y)$ for all $x,y\in Q$.
\item $k(y) = d + g(y)$ for all $x\in Q$, where $d = k(0) = b\cdot ba$.
\item $g$ is a left pseudoautomorphism with left companion $d$ of $(Q,+)$.
\end{enumerate}
\end{lemma}

\begin{proof}
(1) can be proved similarly to \lemref{kepka-5.2-B}(1) using $(F_l)$,
and (2) and (3) follow similarly.
\end{proof}

\begin{lemma}
\lemlabel{kepka-5.5}
\begin{enumerate}
\item $x + (y + z) = (f(x) + y) + (p(x) + z)$ for all $x,y,z\in Q$.
\item $x = f(x) + p(x)$ for all $x\in Q$.
\item $(x + pf\iv(x)) + (y + z) = (x + y) + (pf\iv(x) + z)$ for all $x,y,z\in Q$.
\end{enumerate}
\end{lemma}

\begin{proof}
For (1): We have
\begin{align*}
h(u) + (kh(v) + gk(w)) &= h(u) + k(h(v) + k(w))\\
&= u\cdot vw = uv\cdot \alpha(u)w \\
&= h(h(u) + k(v)) + k(h\alpha(u) + k(w)) \\
&= (fh(u) + hk(v)) + (hk\alpha(u) + gk(w))
\end{align*}
for all $u,v,w\in Q$ by Lemma \lemref{kepka-5.2-A}. Consequently, 
$x + (y + z) = (f(x) + y) + (p(x) + z)$ for all $x,y,z\in Q$.

(2) follows from (1) by taking $y = z = 0$, and (3) follows from combining
(1) and (2).
\end{proof}

\begin{lemma}
\lemlabel{kepka-5.6}
\begin{enumerate}
\item $(x + y) + z = (x + q(z)) + (y + g(z))$ for all $x,y,z\in Q$.
\item $z = q(z) + g(z)$ for all $z\in Q$.
\item $(x + y) + (qg\iv(z) + z) = (x + qg\iv(z)) + (y + z)$ for all $x,y,z\in Q$.
\end{enumerate}
\end{lemma}

\begin{proof}
This is dual to Lemma \lemref{kepka-5.5}.
\end{proof}

\begin{theorem}
\thmlabel{NK-at-last}
$(Q,+)$ is an NK-loop.
\end{theorem}

\begin{proof}
This follows from Theorem \thmref{NK-char} and Lemma \lemref{kepka-5.5}(3)
(or \lemref{kepka-5.6}(3)).
\end{proof}

\begin{lemma}
\lemlabel{kepka-5.9}
$f,g\in \Aut(Q,+)$ and $c,d\in N(Q,+)$.
\end{lemma}

\begin{proof}
By Lemma \lemref{kepka-5.2-A}, $f$ is a right pseudoautomorphism
with companion $c$. By Theorem \thmref{NK-at-last}, $(Q,+)$ is
an NK-loop, so we may apply Lemma \lemref{kepka-3.12}. The rest
is dual.
\end{proof}

\begin{lemma}
\lemlabel{kepka-5.7}
\begin{enumerate}
\item $p(x), q(x)\in K(Q,+)$ for all $x\in Q$.
\item $p(x) - f(x), q(x) - g(x)\in N(Q,+)$ for all $x\in Q$.
\end{enumerate}
\end{lemma}

\begin{proof}
By Theorem \thmref{NK-char}, $pf\iv(x)\in K(Q,+)$ and 
$-x + pf\iv(x)\in N(Q,+)$ for every $x\in Q$. Since $f$
is a permutation, $p(x)\in K(Q,+)$ and so
$p(x) - f(x) = -f(x) + p(x)\in N(Q,+)$. The rest is dual.
\end{proof}

\begin{lemma}
\lemlabel{kepka-5.10}
\begin{enumerate}
\item $-x + f(x), -x + g(x)\in K(Q,+)$ for all $x\in Q$.
\item $x + f(x), x + g(x)\in N(Q,+)$ for all $x\in Q$.
\end{enumerate}
\end{lemma}

\begin{proof}
By Lemmas \lemref{kepka-5.5} and \lemref{kepka-5.7},
$f(x) - x = (x - f(x))\iv = p(x)\iv \in K(Q,+)$ and so 
$-x + f(x) = -x + (f(x) - x) + x = f(x) - x\in K(Q,+)$.
By Lemma \lemref{kepka-5.7}, $x - 2f(x) = p(x) - f(x)\in N(Q,+)$.
Since $(Q,+)$ is an NK-loop, $3f(x)\in N(Q,+)$ by Corollary \cororef{NK-struct}(2).
Thus $x + f(x) = (x - 2f(x)) + 3f(x)\in N(Q,+)$. The remaining assertions
are dual.
\end{proof}

\begin{lemma}
\lemlabel{kepka-5.11}
$hk\alpha(Q)\subseteq K(Q,+)$ and $hk\beta(Q)\subseteq K(Q,+)$.
\end{lemma}

\begin{proof}
This follows from $hk\alpha f\iv = p$, $hk\beta g\iv = q$, and 
Lemma \lemref{kepka-5.7}.
\end{proof}

\begin{lemma}
\lemlabel{kepka-5.13-A}
For all $x\in Q$, $f(d) + fg(x) + c = d + gf(x) + g(c)$.
In particular, $f(d) + c = d + g(c)$.
\end{lemma}

\begin{proof}
Implicitly using Lemma \lemref{kepka-5.9}, we compute
\begin{align*}
f(d) + fg(x) + c &= fk(0) + fg(x) + h(0) = f(k(0) + g(x)) + h(0) = h(k(0) + g(x)) \\
&= hk(x) = kh(x) \\
&= k(f(x) + h(0)) = k(0) + g(f(x) + h(0)) = k(0) + gf(x) + gh(0)\\
&= d + gf(x) + g(c)
\end{align*}
by repeated use of Lemmas \lemref{kepka-5.2-A}, \lemref{kepka-5.2-B}, and \lemref{kepka-5.9}.
The rest follows from taking $x = 0$.
\end{proof}

\begin{lemma}
\lemlabel{kepka-5.13}
\begin{enumerate}
\item If $c\in K(Q,+)$, then $c\in Z(Q,+)$ and $fg = gf$.
\item If $d\in K(Q,+)$, then $d\in Z(Q,+)$ and $fg = gf$.
\end{enumerate}
\end{lemma}

\begin{proof}
The first assertion of (1) is clear from Lemma \lemref{kepka-5.9},
and the rest follows from Lemma \lemref{kepka-5.13-A}:
$f(d) + c + fg(x) = d + g(c) + gf(x)$. (2) is dual to (1). 
\end{proof}

\begin{lemma}
\lemlabel{kepka-5.14}
\begin{enumerate}
\item If $a\in \alpha(Q)$, then $c\in Z(Q,+)$.
\item If $b\in \beta(Q)$, then $d\in Z(Q,+)$.
\end{enumerate}
\end{lemma}

\begin{proof}
For (1), we have $c = ba\cdot a = hk(a)\in hk\alpha(Q)\subset K(Q,+)$,
by Lemma \lemref{kepka-5.11}. Putting this together with Lemma \lemref{kepka-5.9}
gives the desired result. (2) is dual to (1).
\end{proof}

Putting all this together, we have the following.

\begin{theorem}
\thmlabel{kepka-5.15}
Assume that $a\in \alpha(Q)$, $b\in \beta(Q)$ and $\alpha(b) = \beta(a)$
(see Lemma \lemref{kepka-2.3}). Then:
\begin{enumerate}
\item \quad $(Q,+)$ is an NK-loop.
\item \quad $f,g\in \Aut(Q,+)$ and $fg = gf$.
\item \quad $-x + f(x), -x + g(x)\in K(Q,+)$ for every $x\in Q$.
\item \quad $x + f(x), x + g(x)\in N(Q,+)$ for every $x\in Q$.
\item \quad $e = c + d = ba\cdot a + b\cdot ba = ba\cdot ba \in Z(Q,+)$.
\item \quad $xy = f(x) + g(y) + e$ for all $x,y\in Q$.
\end{enumerate}
\end{theorem}

\begin{proof}
Combine Theorem \thmref{NK-at-last} and Lemmas \lemref{kepka-5.9},
\lemref{kepka-5.10}, \lemref{kepka-5.13}, and \lemref{kepka-5.14}.
\end{proof}

\begin{remark}
\remlabel{kepka-5.16}
For $r\in Q$, put $a = \alpha(r)$ and $b = \beta(r)$. Then
$\alpha(b) = \alpha\beta(r) = \beta\alpha(r) = \beta(a)$
by Lemma \lemref{kepka-2.3}, and so the hypotheses of Theorem
\thmref{kepka-5.15} are satisfied in this case. Note that
$0 = ba = \beta(r)\alpha(r)$ and $e = \beta(rr)\alpha(rr)$.
\end{remark}

\begin{lemma}
\lemlabel{kepka-5.17}
\begin{enumerate}
\item If $a,b\in \alpha(Q)$, then $\alpha(Q)\subset K(Q,+)$.
\item If $a,b\in \beta(Q)$, then $\beta(Q)\subset K(Q,+)$.
\item If $a,b\in \alpha\beta(Q) (= \beta\alpha(Q))$, then $\alpha(Q)\cup \beta(Q)\subset K(Q,+)$.
\end{enumerate}
\end{lemma}

\begin{proof}
For (1): we have $a = \alpha(r)$, $b = \beta(s)$, $r,s\in Q$, and
$\alpha(sx\cdot r) = b\alpha(x)\cdot a = hk\alpha(x)\in K(Q,+)$ for
every $x\in Q$ by Lemma \lemref{kepka-5.11}. Since $R_r L_s$ is
a permutation of $Q$, $\alpha(Q)\subseteq K(Q,+)$. 

Now (2) is dual to (1), and (3) follows from combining (1) and (2).
\end{proof}


\section{Quasigroups Linear Over NK-Loops}
\seclabel{linear2}

In this section, let $Q$ be an NK-loop. We denote the underlying sets
of $N(Q,+)$ and $K(Q,+)$ by just $N$ and $K$, respectively.

Let $e\in N$ and let $f,g$ be commuting automorphisms of $(Q,+)$ such
that $-x + f(x),-x + g(x)\in K$ and $x + f(x),x + g(x)\in N$ for all
$x\in Q$. Now define a multiplication on $Q$ by
\[
xy = f(x) + e + g(y)
\]
for $x,y\in Q$. Denote the corresponding quasigroup by $(Q,\cdot)$.

\begin{proposition}
\proplabel{kepka-6.2}
\begin{enumerate}
\item $(Q,\cdot)$ is an F-quasigroup.
\item $\alpha(x) = -g\iv(e) -g\iv f(x) + g\iv(x)$ for every $x\in Q$.
\item $\beta(x) = f\iv(x) - f\iv g(x) - f\iv(e)$ for every $x\in Q$.
\end{enumerate}
\end{proposition}

\begin{proof}
First, $x = x\alpha(x) = f(x) + e + g\alpha(x)$, and so 
$\alpha(x) = -g\iv (e) - g\iv f(x) + g\iv (x)$. Further, by 
Lemma \lemref{kepka-3.14}, $u + (v + w) = (f(u) + v) + ((-f(u) + u) + w)$
for all $u,v,w\in Q$. Setting $u = f(x) + e$, $v = fg(y) = gf(y)$, and
$w = g(e) + g^2(z)$, we get
\begin{align*}
x\cdot yz &= (f(x) + e) + (fg(y) + (g(e) + g^2(z))) \\
&= u + (v + w) = (f(u) + v) + ((-f(u) + u) + w) \\
&= (f^2(x) + f(e) + fg(y)) + ((-f(e) - f^2(x) + f(x) + e) +  (g(e) + g^2(z))) \\
&= f(f(x) + e + g(y)) + ((e - f(e) - f^2(x) + f(x)) + (g(e) + g^2(z)) \\
& = f(xy) + e + g(f(-g\iv(e) - g\iv f(x) + g\iv(x)) + e + g(z))\\
&= f(xy) + e + g(f\alpha(x) + e + g(z)) \\
&= f(xy) + e + g(\alpha(x)z) \\
&= xy\cdot \alpha(x)z
\end{align*}
(Here we have used $-f(e) - f^2(x) + f(x) + e = e - f(e) - f^2(x) + f(x)$,
since $-f^2(x) + f(x)\in K$ and $e - f(e) = -f(e) + e\in K$.) Thus we
have verified $(F_l)$, and the proof of $(F_r)$ is dual to this.
\end{proof}

\begin{lemma}
\lemlabel{kepka-6.3}
Let $(P,+)$ be a subloop of $(Q,+)$ such that $e\in P$ and $f(P) = P = g(P)$.
Then
\begin{enumerate}
\item $(P,\cdot)$ is a subquasigroup of $(Q,\cdot)$.
\item If $(P,+)$ is a normal subloop, then $(P,\cdot)$ is a normal subquasigroup
(and then the corresponding normal congruences coincide).
\end{enumerate} 
\end{lemma}

\begin{proof}
(1) is clear. Now assume that $(P,+)$ is normal in $(Q,+)$ and denote by $\rho$
the corresponding normal congruence of $(Q,+)$; $P$ is a block of $\rho$.
If $(a,b)\in \rho$, then $a-b\in P$, and so $f(a) - f(b) = f(a-b)\in P$,
and so $(f(a),f(b))\in \rho$. Consequently, $(ax,bx) = (f(a) + e + g(x),f(b) + e + g(x))\in \rho$
for each $x\in Q$. The other cases to check are similar, and it follows that 
$\rho$ is a normal congruence of the quasigroup $Q$, too.
\end{proof}

\begin{remark}
\remlabel{kepka-6.4}
Consider the situation from Lemma \lemref{kepka-6.3} and put $(R,+) = (Q,+)/(P,+)
= (Q,+)/\rho$. Then $(R,+)$ is an NK-loop and the automorphisms $f,g$ induce automorphisms
$\bar{f}, \bar{g}$ of $(R,+)$ such that $\pi f = \bar{f}\pi$, $\pi g = \bar{g}\pi$,
where $\pi : Q\to R$ is the natural projection. Moreover, $\bar{e} = \pi(e)\in N(R,+)$.
On the other hand, $(R,\cdot) = (Q,\cdot)/\rho$ is a factor quasigroup of $(Q,\cdot)$
and $\bar{x}\bar{y} = \bar{f}(\bar{x}) + \bar{e} + \bar{g}(\bar{y})$ for all 
$\bar{x},\bar{y}\in R$.
\end{remark}

\begin{lemma}
\lemlabel{kepka-6.5}
Let $(P,+)$ be a subloop of $(Q,+)$ such that either $K\subseteq P$ or
$N\subseteq P$. Then $f(P) = P = g(P)$.
\end{lemma}

\begin{proof}
We have $-x + f(x), -x + f\iv(x), -x + g(x), -x + g\iv(x)\in K$ and
$x + f(x), x + f\iv(x), x + g(x), x + g\iv(x) \in N$.
\end{proof}

\emph{Throughout the rest of this section}, we assume that $e \in Z(Q,+)$.

\begin{lemma}
\lemlabel{kepka-6.7}
\begin{enumerate}
\item $M = M(Q,\cdot) = K (= M(Q,+)$ by Proposition \propref{Mfg_basic}(3)).
\item $\alpha(Q)\cup \beta(Q)\subseteq M$.
\item $x\alpha(z)\cdot yx = xy\cdot \alpha(z)x$ and $x\beta(y)\cdot zx = xz\cdot \beta(y)x$
for all $x,y,z\in Q$.
\item $M$ is a normal subquasigroup of $(Q,\cdot)$.
\item For each $x\in Q$, the subquasigroup $\sbl{x,M}$ is trimedial.
\item $M$ is trimedial.
\item $Q/M$ is a group and in fact, $Q/M \cong N(Q,+)/K(Q,+)$.
\end{enumerate}
\end{lemma}

\begin{proof}
For (1): If $a\in M$, then
\[
(f^2(x) + fg(a)) + (fg(y) + g^2(x)) + c = xa\cdot yx = xy\cdot ax
= (f^2(x) + fg(y)) + (fg(a) + g^2(x)) + c ,
\]
where $c = f(e) + g(e) + e$. Setting $x=0$, we get $fg(a)\in K$.
Since $K$ is characteristic, $a\in K$. We have thus shown $M\subseteq K$.
Similarly, using Lemma \lemref{kepka-3.15}, we may show the other inclusion.

For (2) and (3): We have $-f(x) + x\in K$, and so $\alpha(x) = g\iv(-e-f(x) + x)
\in K$. Similarly, $\beta(x)\in K$.

For (4): Since $K(Q,+)$ is a normal subloop of $(Q,+)$ (by Corollary \cororef{NK-struct}(3)),
$(M,\cdot)$ is a normal subquasigroup by Lemmas \lemref{kepka-6.5} and \lemref{kepka-6.3}.

For (5): Let $(P,+)$ be the subloop generated by $\{x,M\}$. By Lemma \lemref{kepka-3.13}(1),
$(P,+)$ is a commutative loop and by Lemmas \lemref{kepka-6.5} and \lemref{kepka-6.3},
$(P,\cdot)$ is a subquasigroup of $(Q,\cdot)$. Now $(P,\cdot)$ is trimedial by 
Proposition \propref{kepka-2.8}.

(6) follows from (5).

For (7): From (2), $(Q/M,\cdot)$ is a loop, and hence a group by Lemma \lemref{kepka-2.5}.
Now consider the situation from Remark \remref{kepka-6.4} where $P = M$. Since 
$-x + f(x)\in M$, we have $\bar{f} = \Id_{Q/M}$. Similarly, $\bar{g} = \Id_{Q/M}$
and, since $e\in M$, have $\bar{e} = \bar{0}$ and
$\bar{x}\bar{y} = \bar{x} + \bar{y}$ for all $\bar{x},\bar{y}\in R = Q/M$. Thus
$R = (Q,+)/K(Q,+) \cong N(Q,+)/(K(Q,+)\cap N(Q,+)) = N(Q,+)/Z(Q,+)$.
\end{proof}

\begin{corollary}
\corolabel{kepka-6.8}
$(Q,\cdot)$ is monomedial.
\end{corollary}

\begin{lemma}
\lemlabel{kepka-6.9}
\begin{enumerate}
\item $(N,\cdot)$ is a normal subquasigroup of $(Q,\cdot)$.
\item $(N,\cdot)$ is an FG-quasigroup.
\item For all $x,y\in Q$, the subquasigroup generated by $\{x,y\}\cup N$
is an FG-quasigroup.
\item $(Q/N,\cdot)$ is a symmetric, distributive quasigroup (in particular,
every block of the congruence corresponding to $N$ is a subquasigroup of $(Q,\cdot)$).
\end{enumerate}
\end{lemma}

\begin{proof}
For (1): Combine Lemmas \lemref{kepka-6.3} and \lemref{kepka-6.5}, and
Corollary \cororef{NK-struct}(2).

For (2) and (3): The subloop generated by the set is a group.

For (4): Consider again the situation from Remark \remref{kepka-6.4} where
$P = N(Q,+)$. Since $(Q,+)$ is an NK-loop, the factor loop $(R,+)
= (Q,+)/(N,+)$ is a commutative Moufang loop (Corollary \cororef{NK-struct}).
Further, $x + f(x)\in N$ so that $\bar{f}(\bar{x}) = -\bar{x}$. Similarly,
$\bar{g}(\bar{x}) = -\bar{x}$ and $\bar{e}=\bar{0}$. Thus
$\bar{x}\bar{y} = -\bar{x}-\bar{y}$ for all $\bar{x},\bar{y}\in R$ and
we apply Proposition \propref{kepka-2.11}. Finally, every block of $\rho$ 
is a subquasigroup, since $R = Q/\rho = Q/N$ is idempotent.
\end{proof}

\begin{lemma}
\lemlabel{kepka-6.10}
Let $a,b,c,d\in Q$. Then $ab\cdot cd = ac\cdot bd$ if and only if
$(a + b) + (c + d) = (a + c) + (b + d)$.
\end{lemma}

\begin{proof}
It is easy to see that $ab\cdot cd = ac\cdot bd$ if and only if
$(fg\iv(a) + b) + (c + gf\iv(d)) = (fg\iv(a) + c) + (b + gf\iv(d))$.
Since $a - fg\iv(a)\in N$ and $-gf\iv(d) + d\in N$, the latter equality
is equivalent to $(a + b) + (c + d) = (a + c) + (b + d)$.
\end{proof}

\begin{lemma}
\lemlabel{kepka-6.11}
If $a\in N$ and $b\in M$, then $xa\cdot by = xb\cdot ay$ for all $x,y\in Q$.
\end{lemma}

\begin{proof}
The equality follows easily from Lemma \lemref{kepka-6.10}.
\end{proof}

\begin{lemma}
\lemlabel{kepka-6.12}
The quasigroup $(Q,\cdot)$ is a homomorphic image of the direct product
$N\times M$ of the quasigroups $N$ and $M$.
\end{lemma}

\begin{proof}
According to Lemma \lemref{kepka-6.11}, the mapping $(a,b)\mapsto ab$ is
a homomorphism of $N\times M$ into $Q$. On the other hand, if $x\in Q$,
then $x = c + d$, $c\in N$, $d\in M$, and $x = f\iv(c)\cdot g\iv(d - e)$,
where $f\iv(c)\in N$, $g\iv(d-e)\in M$. Thus $Q = NM$ and the homomorphism
is a projection.
\end{proof}

\begin{lemma}
\lemlabel{kepka-6.13}
Every (at most) three-generated subquasigroup of $(Q,\cdot)$ is an FG-quasi\-group.
\end{lemma}

\begin{proof}
F-quasigroups with the indicated property form an equational class of quasigroups,
and this class contains all FG-quasigroups and all trimedial quasigroups. Our
result now follows from Lemma \lemref{kepka-6.12}.
\end{proof}

\begin{lemma}
\lemlabel{kepka-6.14}
\begin{enumerate}
\item $Z = N\cap M$ is a normal subquasigroup of $Q$.
\item $(Z,\cdot)$ is a medial quasigroup.
\item For every $x\in Q$, the subquasigroup generated by the set $\{x\}\cup Z$
is medial.
\item $Q/Z$ is isomorphic to a subquasigroup of $Q/N\times Q/M$, which is the
product of a group and a symmetric distributive quasigroup.
\end{enumerate}
\end{lemma}

\begin{proof}
This follows from Lemmas \lemref{kepka-6.7} and \lemref{kepka-6.9}.
\end{proof}

\begin{lemma}
\lemlabel{kepka-6.15}
The following conditions are equivalent.
\begin{enumerate}
\item $(Q,+)$ is commutative.
\item $(Q,\cdot)$ is trimedial.
\item $(Q,\cdot)$ is dimedial.
\item $xx\cdot yx = xy\cdot xx$ for all $x,y\in Q$.
\item $xx\cdot yy = xy\cdot xy$ for all $x,y\in Q$.
\item $N(Q,+)$ is an abelian group.
\end{enumerate}
\end{lemma}

\begin{proof}
(1)$\implies$(2): This follows from Proposition \propref{kepka-2.8}.

(2)$\implies$(3) and (3)$\implies$(4),(5): trivial.

(4)$\implies$(1): By Lemma \lemref{kepka-6.10}, $(x + x) + (y + x)
= (x + y) + (x + x)$ for all $x,y\in Q$. By diassociativity, we 
may cancel to get $x + y = y + x$, i.e., $(Q,+)$ is commutative.

(5)$\implies$(1): This is proved similarly to the previous case.

(1)$\implies$(6): trivial.

(6)$\implies$(2): $(N,\cdot)$ is a trimedial quasigroup, and so is
$N\times M$, and so $Q$ is trimedial by Lemma \lemref{kepka-6.12}.
\end{proof}

\begin{lemma}
\lemlabel{kepka-6.16}
The following conditions are equivalent.
\begin{enumerate}
\item $(Q,+)$ is a group.
\item $(Q,\cdot)$ is an FG-quasigroup.
\item Every four-generated subquasigroup of $Q$ is an FG-quasigroup.
\item $x\alpha(u)\cdot yz = xy \cdot \alpha(u)z$ for all $x,y,z,u\in Q$.
\item $x\beta(u)\cdot yz = xy \cdot \beta(u)z$ for all $x,y,z,u\in Q$.
\item $x\alpha(u)\cdot \beta(v)z = x\beta(v) \cdot \alpha(u)z$ for all $x,y,u,v\in Q$.
\item $K(Q,+)$ is an abelian group.
\end{enumerate}
\end{lemma}

\begin{proof}
(1)$\iff$(2), (2)$\implies$(3), (4)$\implies$(6), and (5)$\implies$(6)
are all trivial.

(3)$\implies$(4): Four letters occur in the equality in (4), and so we may assume
without loss of generality that (2), and hence (1), hold. In view of Lemma
\lemref{kepka-6.10}, we have to show that $x + \alpha(u) + y + z =
x + y + \alpha(u) + z$, i.e., $\alpha(u) + y = y + \alpha(u)$. However, by 
Lemma \lemref{kepka-6.7}(2), $\alpha(u)\in K(Q,+) = Z(Q,+)$. 

(6)$\implies$(7): By Lemma \lemref{kepka-6.10}, $(x + \alpha(u)) + (\beta(v) + y)
= (x + \beta(v)) + (\alpha(u) + y)$ for all $x,y,u,v\in K$. Setting $v =0$ and
taking into account that $e\in Z(K) = Z(K(Q,+))$, we get $x + (-g\iv f(u) + g\iv(u)) + y
= x + ((-g\iv f(u) + g\iv(u)) + y)$, i.e., $-g\iv f(u) + g\iv(u)\in Z(K)$ for
every $u\in K$. Then $-f(u) + u\in Z(K)$ and since $u + f(u)\in Z(K)$,
we get $2u\in Z(K)$. But $3u\in Z(K)$ implies $u\in Z(K)$. Thus $(K,+)$ is an
abelian group.

(7)$\implies$(1): By Corollary \cororef{NK-struct}(1), $(Q,+)$ is an image of
the product $N(Q,+)\times K(Q,+)$. This product is a group.

\end{proof}

\begin{lemma}
\lemlabel{kepka-6.17}
The following conditions are equivalent.
\begin{enumerate}
\item $(Q,+)$ is an abelian group.
\item $(Q,\cdot)$ is medial.
\item $xx\cdot yx = xy\cdot xx$ and $x\alpha(u)\cdot \beta(v)y = x\beta(v)\cdot \alpha(u)y$
for all $x,y,u,v\in Q$.
\item $xx\cdot yy = xy\cdot xy$ and $x\alpha(u)\cdot \beta(v)y = x\beta(v)\cdot \alpha(u)y$
for all $x,y,u,v\in Q$.
\end{enumerate}
\end{lemma}

\begin{proof}
Combine Lemmas \lemref{kepka-6.15} and \lemref{kepka-6.16}.
\end{proof}

\begin{lemma}
\lemlabel{kepka-6.18}
If $(Q,\cdot)$ is unipotent, then $(Q,\cdot)$ is medial.
\end{lemma}

\begin{proof}
According to Lemma \lemref{kepka-6.17}, it is sufficient to show that $(Q,+)$
is an abelian group. We have $f(x) + g(x) + e = xx = 0\cdot 0 = e$ and
so $f(x) + g(x) = x$ for every $x\in Q$. Then $f = -g$, and so the mapping
$x\mapsto -x$ is an automorphism of $(Q,+)$. Hence $(Q,+)$ is commutative.
Further, from Lemma \lemref{kepka-6.9}(4), $Q/N$ is both unipotent and
idempotent. Then $Q/N$ is trivial, so that $Q = N$ and $(Q,+)$ is an abelian group.
\end{proof}


\section{Quasigroups Linear Over NK-Loops II}
\seclabel{linear3}

We continue with the notational conventions of the preceding section, and
continue to assume that $e\in Z(Q,+)$.

\begin{lemma}
\lemlabel{kepka-7.2}
Let $p,q : Q\to Q$ be mappings. Then:
\begin{enumerate}
\item $(p,q)\in \AAA(Q)$ if and only if there exists $r\in N(Q,+)$ such
that $p(x) = f(r) + x$ and $q(x) = r + x$ for every $x\in Q$.
\item $(p,q)\in \BBB(Q)$ if and only if there exists $r\in N(Q,+)$ such
that $p(x) = x + g(r)$ and $q(x) = x + r$ for every $x\in Q$.
\item $(p,q)\in \CCC(Q)$ if and only if there exists $r\in N(Q,+)$ such
that $p(x) = x + r$ and $q(x) = g\iv f(r) + x$ for every $x\in Q$.
\end{enumerate}
\end{lemma}

\begin{proof}
For (1): $(p,q)\in \AAA(Q)$ if and only if $p(f(x) + g(y) + e)
= p(xy) = q(x)y = fq(x) + g(y) + e$, or equivalently, if and
only if $p(x+y) = fqf\iv(x) + y$ for all $x,y\in Q$, that is,
if and only if $(p,fqf\iv)\in \AAA(Q,+)$. The rest is easy.

The proofs of (2) and (3) are similar.
\end{proof}

\begin{lemma}
\lemlabel{kepka-7.3}
\begin{enumerate}
\item $\AAA_l(Q) = \AAA_r(Q) = \CCC_l(Q)$.
\item $\BBB_l(Q) = \BBB_r(Q) = \CCC_r(Q)$.
\item The permutation groups $\AAA_l(Q)$ and $\BBB_l(Q)$
are isomorphic to the group $N(Q,+)$.
\end{enumerate}
\end{lemma}

\begin{proof}
This follows from definitions and Lemma \lemref{kepka-7.2}.
\end{proof}

\begin{lemma}
\lemlabel{kepka-7.4}
Let $\rho$ be the normal congruence of $(Q,\cdot)$ (and $(Q,+)$ as well)
corresponding to $N$ (see Lemma \lemref{kepka-6.9}). Then:
\begin{enumerate}
\item $(a,b)\in \rho$ if and only if $a = p(b)$ for some $p\in \AAA_l(Q)$.
\item $(a,b)\in \rho$ if and only if $a = p(b)$ for some $p\in \BBB_l(Q)$.
\item $Q/\rho = \{ N + u : u\in K \} = \{ Nu : u \in K\} = \{ uN : u\in K\}$.
\end{enumerate}
\end{lemma}

\begin{proof}
This is elementary using Lemmas \lemref{kepka-7.2}, \lemref{kepka-7.3}, and
\lemref{kepka-6.9}.
\end{proof}

\begin{construction}
\conslabel{kepka-7.5}
Fix $a\in K$, $b\in N$, and set $\tau(x) = (x + b) + a$
($= x + b + a = x + a + b = (x + a) + b$) and 
$x\ast y = ((x - b) + y) - a$ ($= (x - b + y) - a =
-a + (x - b + y))$ for all $x,y\in Q$. We have defined
a new binary operation $\ast : Q\times Q\to Q$. 

\noindent (i) Using $a\in K$, we get
\begin{align*}
\tau(x)\ast \tau(y) &= ((x+b+a)-b)+(y+b+a))-a = 
((x+b) +(y+b+a))-a \\
&= ((x + y + b) + 2a) - a = (x + y + b) - a = \tau(x + y)
\end{align*}
for $x,y\in Q$. This $\tau : (Q,+)\to (Q,\ast)$ is an isomorphism
of binary structures. In particular, $(Q,\ast)$ is an NK-loop and
the neutral element of $(Q,\ast)$ is $a + b = b + a$.

\noindent (ii) Since $\tau$ is an isomorphism, we have
$N(Q,\ast) = \tau(N(Q,+)) = N(Q,+) + a$ and
$K(Q,\ast) = \tau(K(Q,+)) = K(Q,+) + b$. 

\noindent (iii) The mapping $h = \tau f \tau\iv$ is an
automorphism of $(Q,\ast)$; we have
\[
h(x) = \tau f((x-b)-a) = \tau(f(x)-f(b)-f(a))
= (f(x)-f(b)-f(a))+b+a = f(x) + (b-f(b)) + (a-f(a))
\]
using $b-f(b), a+f(a)\in Z(Q,+)$. Similarly, $k = \tau g \tau\iv$
is an automorphism of $(Q,\ast)$, and we have
$k(x) = g(x) + (b-g(b)) + (a-g(a))$.

\noindent (iv) Now $hk\tau = h\tau g = \tau fg = \tau g f
= k \tau f = kh\tau$, and so $hk = kh$.

\noindent (v) For each $x\in Q$, $\tau\iv (x\ast h(x))
= \tau\iv(x) + \tau\iv h(x) = \tau\iv(x) + f\tau\iv(x)\in N(Q,+)$,
and so $x * h(x)\in \tau(N(Q,+)) = N(Q,\ast)$. Similarly,
$x \ast k(x)\in N(Q,\ast)$.

\noindent (vi) If $\tilde{x}$ denotes the inverse of $x$ in the
loop $(Q,\ast)$, then $\tilde{x}\ast h(x)\in K(Q,\ast)$ and
$\tilde{x}\ast k(x)\in K(Q,\ast)$. 

\noindent (vii) We have
$xy = f(x) + g(y) + e = (h(x) + (f(a) - a) + (f(b) - b)) 
+ (k(y) + (g(a) - a) + (g(b) - b)) + e =
(h(x) + (f(a)-a)) + (k(y) + (g(a)-a)) + r$
where $r = (f(b)-b) + (g(b)-b) + e\in Z(Q,+)$.
Henceforth, $xy = ((h(x) + k(y)) - a) + s$ where
$s = 3a + (a + f(a)) + (a + g(a)) + r\in Z(Q,+)$.
Finally, 
$e_1 = 2b + a + s \in N(Q,+) + a = N(Q,\ast)$ and
$h(x)\ast e_1 \ast k(y) = 
(((((h(x) - b) + e_1) - a) - b) + k(y)) - a =
((h(x) + (-b + e_1 - a- b)) + k(y)) - a =
((h(x) + k(y)) - a) + s$. We have shown that
$xy = h(x) \ast e_1 \ast k(y)$ for all $x,y\in Q$.

\noindent (viii) Note that $e_1 \in Z(Q,\ast)$ if and
only if $b\in Z(Q,+)$ (or $b\in K(Q,+)$).

\noindent (ix) Put $N_1 = N + a$. By Lemma \lemref{kepka-6.9},
$N_1$ is a block of the congruence $\rho$ corresponding to $N$
and $N_1$ is a normal subquasigroup of $(Q,\cdot)$. Now
$N_1$ is the underlying set of $N(Q,\ast)$ and, by (vii),
$uv = h(u)\ast e_1 \ast k(v)$ for all $u,v\in N_1$. In
particular, $N_1$ is an FG-quasigroup isotopic to $N(Q,\ast)$.
The groups $N(Q,+)$ and $N(Q,\ast)$ are isomorphic.

\noindent (x) Put $e_2 = b + s$, so that $\tau(e_2) = e_1$
and $e_2\in N(Q,+)$. Now define a binary operation $\Delta :Q\times Q\to Q$
by $x\Delta y = f(x) + e_2 + g(y)$ for $x,y\in Q$.
Then $\tau(x\Delta y) = f\tau(x) \ast \tau(e_2)\ast g\tau(y)
= h\tau(x) \ast e_1 \ast k\tau(y) = \tau(x)\tau(y)$.
Thus $\tau : (Q,\Delta)\to (Q,\cdot)$ is an isomorphism of
quasigroups.

\noindent (xi) Put $e_3 = (e + b) + a$ so that $\tau(e) = e_3$
and $e_3\in N(Q,\ast)$. Now define a binary operation
$\nabla : Q\times Q\to Q$ by $x\nabla y = h(x) \ast e_3\ast k(y)$.
Then $\tau(xy) = \tau(x)\nabla \tau(y)$ and consequently,
$\tau : (Q,\cdot)\to (Q,\nabla)$ is an isomorphism of quasigroups.
\end{construction}


\section{Arithmetic Forms of F-Quasigroups}
\seclabel{forms}

An ordered five-tuple $(Q,+,f,g,e)$ will be called an \emph{arithmetic
form} of a quasigroup $(Q,\cdot)$ if

\noindent (1) \qquad $(Q,+)$ is an NK-loop;

\noindent (2) \qquad $f,g$ are commuting automorphisms of $(Q,+)$;

\noindent (3) \qquad $x+f(x),x +g(x)\in N(Q,+)$ for every $x\in Q$;

\noindent (4) \qquad $-x+f(x),-x+g(x)\in K(Q,+)$ for every $x\in Q$;

\noindent (5) \qquad $e\in N(Q+)$;

\noindent (6) \qquad $xy = f(x) + e + g(y)$ for all $x,y\in Q$.

If, moreover,

\noindent (7) \qquad $e\in Z(Q,+)$,

\noindent then the arithmetic form will be called \emph{strong}.

\begin{theorem}
\thmlabel{kepka-8.2}
The following conditions are equivalent for a quasigroup $(Q,\cdot)$.
\begin{enumerate}
\item $(Q,\cdot)$ is an F-quasigroup.
\item $(Q,\cdot)$ has at least one strong arithmetic form.
\item $(Q,\cdot)$ has at least one arithmetic form.
\end{enumerate}
\end{theorem}

\begin{proof}
(1)$\implies$(2): Take $r\in Q$ arbitrarily and put $a = \alpha(r)$,
$b = \beta(r)$ (see Remark \remref{kepka-5.16}). Then 
$\alpha(b) = \beta(a)$ and, by Theorem \thmref{kepka-5.15}, we get
a strong arithmetic form $(Q,+,f,g,e)$ of the quasigroup $(Q,\cdot)$.
Note that $0 = ba = \alpha(r)\beta(r)$ and $e = \beta(rr)\alpha(rr)$.

(2)$\implies$(3): trivial

(3)$\implies$(1): This follows from Proposition \propref{kepka-6.2}.
\end{proof}

\begin{lemma}
\lemlabel{kepka-8.3}
Let $(Q,+,f,g,e_1)$ and $(P,+,h,k,e_2)$ be arithmetic forms of F-quasigroups
$(Q,\cdot)$ and $(P,\cdot)$, respectively. Let $\varphi : Q\to P$ be a
mapping such that $\varphi( 0_Q ) = 0_P$. Then $\varphi$ is a homomorphism
of quasigroups if and only if $\varphi$ is a homomorphism of loops such
that $\varphi f = h \varphi$, $\varphi g = k\varphi$ and $\varphi(e_1) = e_2$.
\end{lemma}

\begin{proof}
Assume that $\varphi$ is a homomorphism of quasigroup structures, the
other case being easy. Now $\varphi(f(x) + e_1 + g(y)) = \varphi(xy)
= \varphi(x)\varphi(y) = h\varphi(x) + e_2 + k\varphi(y)$ for all $x,y\in Q$.
Setting $x = y = 0_Q$, we get $\varphi(e_1) = e_2$. Setting $y = g\iv(-e_1)$,
we get $\varphi f(x) = h\varphi(x) + \varphi(e_1) + k\varphi g\iv(-e_1)$,
and hence $x = 0_Q$ yields $\varphi(e_1) + + k\varphi g\iv(-e_1) = 0$.
Thus $\varphi f = h \varphi$, and similarly, $\varphi g = k\varphi$. From
this we conclude that 
$\varphi(x + e_1 + y) = \varphi(x) +\varphi(e_1) + \varphi(y)$ for all
$x,y\in Q$. In particular, $\varphi(e_1 + y) = \varphi(e_1)+ \varphi(y)$,
$\varphi(x + r) = \varphi(x) + \varphi(r)$, $r = e_1 + y$, and so
$\varphi : (Q,+)\to (P,+)$ is a homomorphism.
\end{proof}

\begin{lemma}
\lemlabel{kepka-8.4}
Let $(Q,+,f_1,g_1,e_1)$ and $(Q,\ast,f_2,g_2,e_2)$ be arithmetical forms of
an F-quasigroup $(Q,\cdot)$ such that the loops $(Q,+)$ and $(Q,\ast)$ have
the same neutral element $0$. Then $(Q,+) = (Q,\ast)$, $f_1 = f_2$, $g_1 = g_2$,
$e_1 = e_2$, i.e., the forms coincide.
\end{lemma}

\begin{proof}
The assertion follows from Lemma \lemref{kepka-8.3} where $\varphi = \Id_Q$.
\end{proof}

\begin{theorem}
\thmlabel{kepka-8.5}
Let $(Q,\cdot)$ be an F-quasigroup. Then there exists a one-to-one correspondence
between (strong) arithmetic forms of the quasigroup and elements from $Q$
(resp. $M(Q)$). More precisely, for every element $w\in Q$ ($w\in M(Q)$) there
exists just one arithemetic form of $Q$ such that $w$ is a neutral element of
the corresponding loop.
\end{theorem}

\begin{proof}
By Theorem \thmref{kepka-8.2}(ii), $(Q,\cdot)$ has at least one strong arithmetic
form, say $(Q,+,f,g,e)$. Now if $w\in Q$, then $w = a + b$ for some $a\in K(Q,+)$
and some $b\in N(Q,+)$. By Construction \consref{kepka-7.5}, we get an arithmetic
form $(Q,\ast,h,k,e_1)$ of $Q$ such that $0_{(Q,*)} = a + b = w$. (By \consref{kepka-7.5}(viii)
and Lemma \lemref{kepka-6.7}(1), the form is strong iff $b\in Z(Q,+)$ and hence iff
$w\in M(Q)$.) Finally the uniqueness of the form $(Q,\ast,h,k,e_1)$ follows from
Lemma \lemref{kepka-8.4}.
\end{proof}

Denote by $\mathcal{F}_p$ the equational class (and category)
of pointed $F$-quasigroups. That is, $\mathcal{F}_p$ consists of ordered pairs
$(Q,a)$, where $Q$ is an F-quasigroup and $a \in Q$.
If $(P,b)\in \mathcal{F}_p$, then a mapping $\varphi : Q\to P$ is
a homomorphism in $\mathcal{F}_p$ if and only if $\varphi$ is a 
homomorphism of quasigroups and $\varphi(a) = b$. Further, 
put $\mathcal{F}_m = \{(Q,a) \in \mathcal{F}_p : a \in M(Q)\}$.
Then $\mathcal{F}_m$ is an equational subclass (and a full subcategory)
of $\mathcal{F}_p$. 

Denote by $\mathcal{E}$ the equational class (and category again) of
algebras $(Q,+,f,g,f\iv,g\iv,e)$ where $(Q,+)$ is an NK-loop and
the conditions (2),(3),(4),(5) of the definition of arithmetic form hold.
If $(P,+,h,k,h\iv,k\iv,e_1)\in \mathcal{E}$, then a mapping
$\varphi : Q\to P$ is a homomorphism in $\mathcal{E}$ if and only if
$\varphi$ is a homomorphism of loops such that $\varphi f = h\varphi$,
$\varphi g = k\varphi$, and $\varphi(e) = e_1$. Further, put 
$\mathcal{E}_c = \{ (Q,+,f,g,f\iv,g\iv,e)\in \mathcal{E} : e\in Z(Q,+)\}$.
Then $\mathcal{E}_c$ is an equational subclass (and full subcategory)
of $\mathcal{E}$.

Let $(Q,a)\in \mathcal{F}_p$. By Theorem \thmref{kepka-8.5} (and its proof),
there is just one arithmetic form $(Q,+,f,g,e)$ of $(Q,\cdot)$ such that
$a = 0_{(Q,+)}$ (and $e\in Z(Q,+)$ if and only if $a\in M(Q)$). Now
put $\Phi((Q,a)) = (Q,+,f,g,f\iv,g\iv,e)$.

Let $(Q,+,f,g,f\iv,g\iv,e)\in \mathcal{E}$. Then set
$\Psi((Q,+,f,g,f\iv,g\iv,e)) = (Q, 0_{(Q,+)})\in \mathcal{F}_p$, where
a multiplication on $Q$ is defined by $xy = f(x) + e + g(y)$ 
($0\in M(Q)$ if and only if $e\in Z(Q,+)$).

It follows from Theorem \thmref{kepka-8.5} that $\Phi$ and $\Psi$ are
mutually inverse, one-to-one correspondences between the classes
$\mathcal{F}_p$ and $\mathcal{E}$, $\Phi : \mathcal{F}_p\to \mathcal{E}$
and $\Psi : \mathcal{E}\to \mathcal{F}_p$. If $A,B\in \mathcal{F}_p$,
then a homomorphism $\varphi : A\to B$ is a homomorphism in $\mathcal{F}_p$
if and only if $\varphi : \Phi(A)\to \Phi(B)$ is a homomorphism in
$\mathcal{E}$ (see Lemma \lemref{kepka-8.3}). This implies that the
classes $\mathcal{F}_p$ and $\mathcal{E}$ are equivalent. 

Summarizing this discussion, we have the following.

\begin{corollary}
\corolabel{kepka-8.7}
The class $\mathcal{F}_p$ of pointed F-quasigroups is equivalent to
the class $\mathcal{E}$. The equivalence restricts to an equivalence
between the class $\mathcal{F}_m$ of $M$-pointed F-quasigroups and
the class $\mathcal{E}_c$.
\end{corollary}

\begin{remark}
\remlabel{kepka-8.8}
Let $\varphi : Q\to P$ be a homomorphism of F-quasigroups. If $a\in \alpha(Q)$
($a\in \beta(Q)$, resp.), then $\varphi(a)\in \alpha(P)$ ($\varphi(a)\in \beta(P)$,
resp.), and $(Q,a), (P,\varphi(a))\in \mathcal{F}_m$ and 
$\varphi : (Q,a)\to (P,\varphi(a))$ is a homomorphism in the class $\mathcal{E}_c$.
\end{remark}


\section{Summary of Structure Results on F-Quasigroups}
\seclabel{summary}

\begin{theorem}
\thmlabel{kepka-9.1}
Let $(Q,\cdot)$ be an F-quasigroup.
\begin{enumerate}
\item $M = M(Q)$ is a normal subquasigroup of $Q$.
\item $M$ is a trimedial quasigroup and $Q/M$ is a group.
\item $\alpha(Q)\cup \beta(Q)\subseteq M$.
\item For each $a\in Q$, the subquasigroup generated by the set $\{a\}\cup M$
is trimedial.
\end{enumerate}
\end{theorem}

\begin{proof}
Combine Theorem \thmref{kepka-8.2} with Lemma \lemref{kepka-6.7}.
\end{proof}

\begin{corollary}
\corolabel{kepka-9.2}
Let $(Q,\cdot)$ be an F-quasigroup.
\begin{enumerate}
\item Both $\alpha(Q)$ and $\beta(Q)$ are trimedial subquasigroups. Moreover,
the subquasigroup generated by $\alpha(Q)\cup \beta(Q)$ is trimedial.
\item $x\alpha(z)\cdot yx = xy\cdot \alpha(z)x$ and 
$x\beta(z)\cdot yx = xy\cdot \beta(z)x$ for all $x,y,z\in Q$.
\item $(Q,\cdot)$ is monomedial.
\end{enumerate}
\end{corollary}

\begin{corollary}
\corolabel{kepka-9.3}
Every one-generated F-quasigroup is medial.
\end{corollary}

\begin{theorem}
\thmlabel{kepka-9.4}
Let $(Q,\cdot)$ be an F-quasigroup. Define a relation $\rho$ on $Q$ by
$(a,b)\in \rho$ if and only if $a = p(b)$ for a regular permutation 
$p$ of $Q$. Then:
\begin{enumerate}
\item $\rho$ is a normal congruence of $(Q,\cdot)$.
\item Every block of $\rho$ is a normal subquasigroup of $Q$ and
and an FG-quasigroup.
\item $Q/\rho$ is a symmetric distributive quasigroup.
\item If $N$ is a block of $\rho$ and $a,b\in Q$, then the subquasigroup
generated by the set $\{a,b\}\cup N$ is an FG-quasigroup.
\item Every (at most) three-generated subquasigroup of $(Q,\cdot)$ is
an FG-quasigroup.
\end{enumerate}
\end{theorem}

\begin{proof}
See Theorem \thmref{kepka-8.2}, Lemmas \lemref{kepka-7.4} and \lemref{kepka-6.9},
and Construction \consref{kepka-7.5}(ix).
\end{proof}

\begin{corollary}
\corolabel{kepka-9.5}
Every (at most) three-generated F-quasigroup is an FG-quasigroup.
\end{corollary}

\begin{theorem}
\thmlabel{kepka-9.6}
Let $(Q,\cdot)$ be an F-quasigroup and $N$ a block of the normal congruence $\rho$
(see Theorem \thmref{kepka-9.4}). Then the mapping $(a,u)\mapsto au$, $a\in N$,
$u\in M$ is a surjective homomorphism of the direct product $N\times M(Q)$
onto $Q$.
\end{theorem}

\begin{proof}
See Theorem \thmref{kepka-8.2} and Lemma \lemref{kepka-6.12}.
\end{proof}

\begin{proposition}
\proplabel{kepka-9.7}
The following conditions are equivalent for an F-quasigroup $(Q,\cdot)$.
\begin{enumerate}
\item $(Q,\cdot)$ is trimedial.
\item $(Q,\cdot)$ is dimedial.
\item $xx\cdot yx = xy\cdot xx$ for all $x,y\in Q$.
\item $xx\cdot yy = xy\cdot xy$ for all $x,y\in Q$.
\item At least one of the blocks of the normal congruence $\rho$
is a trimedial quasigroup.
\end{enumerate}
\end{proposition}

\begin{proof}
See Theorem \thmref{kepka-8.2}, Lemma \lemref{kepka-6.15}, and
Theorem \thmref{kepka-9.6}.
\end{proof} 

\begin{proposition}
\proplabel{kepka-9.8}
The following conditions are equivalent for an F-quasigroup $(Q,\cdot)$.
\begin{enumerate}
\item $Q$ is an FG-quasigroup.
\item Every four-generated subquasigroup of $Q$ is an FG-quasigroup.
\item $x\alpha(u)\cdot yz = xy\cdot \alpha(u)z$ for all $x,y,z,u\in Q$.
\item $x\beta(u)\cdot yz = xy\cdot \beta(u)z$ for all $x,y,z,u\in Q$.
\item $x\alpha(u)\cdot \beta(v)y = x\beta(v)\cdot \alpha(u)y$ for all $x,y,u,v\in Q$.
\item $M(Q)$ is a medial quasigroup.
\end{enumerate}
\end{proposition}

\begin{proof}
See Theorem \thmref{kepka-8.2} and Lemma \lemref{kepka-6.16}.
\end{proof} 

\begin{remark}
\remlabel{kepka-9.9}
Let $(Q,\cdot)$ be an F-quasigroup. By Corollary \cororef{kepka-9.2}(1), the
subquasigroup generated by $\alpha(Q)\cup \beta(Q)$ is trimedial. 
In particular, $Q$ is trimedial provided that $\alpha(Q) = Q$ or $\beta(Q) = Q$.
On the other hand, if $\alpha$ ($\beta$, resp.) is injective, then 
$Q$ can be imbedded into an F-quasigroup $(P,\cdot)$ such that
$\alpha$ ($\beta$, resp.) is a permutation of $P$. (This is a standard
construction using the fact that $\alpha$ and $\beta$ are endomorphisms.)
But then $(P,\cdot)$ is trimedial, and hence so is $(Q,\cdot)$.
\end{remark}

\begin{remark}
\remlabel{kepka-9.10}
It follows immediately from Theorem \thmref{kepka-9.1} that every 
(congruence-) simple F-quasigroup is either a simple group or a
simple trimedial quasigroup. While there are many simple groups,
simple trimedial quasigroups are necessarily finite and medial,
and all of them can be found in \cite{JKN}.
\end{remark}

\begin{remark}
\remlabel{kepka-9.11}
Let $(Q,\cdot)$ be an F-quasigroup and $N_1, N_2$ two blocks of the
normal congruence $\rho$ (see Theorem \thmref{kepka-9.4}). Then both
$N_1$ and $N_2$ are FG-quasigroups and it follows from Construction
\consref{kepka-7.5}(ix) that these quasigroups are isotopic to
isomorphic groups.
\end{remark}

\begin{remark}
\remlabel{kepka-9.12}
Let $(Q,+,f,g,e)$ be a strong arithmetic form of an F-quasigroup $(Q,\cdot)$.
If $x\circ y = f(x) + g(y)$ for all $x,y\in Q$, then $(Q,\circ)$ is again
an F-quasigroup and the neutral element $0$ is an idempotent element
of $(Q,\circ)$. We have $x\circ y = (xy) - e$, so that the quasigroups
$(Q,\cdot)$ and $(Q,\circ)$ are isotopic. In this way, we have proved
that every F-quasigroup is isotopic to an F-quasigroup containing at
least one idempotent element. 
\end{remark}

\begin{remark}
\remlabel{kepka-9.13}
\noindent (i) Consider a finite F-quasigroup $(Q,\cdot)$ which is
minimal with respect to the property of 
\emph{not} being an FG-quasigroup. It follows from Theorem
\thmref{kepka-9.6} that $M(Q)$ is not an FG-quasigroup, and hence
$Q = M(Q)$ is trimedial. Of course, $Q$ is not medial.  Now according
to \cite{K-nonmed}, we have $|Q| = 81$. By \cite{K-nonmed}, there exist
just 35 isomorphism classes of nonmedial, trimedial quasigroups of order $81$. 

\smallskip

\noindent (ii) Consider a finite F-quasigroup $(P,\cdot)$ which is
minimal with respect to the property of 
\emph{not} being trimedial. It follows easily from Theorem \thmref{kepka-9.6}
that $P$ is a copy of the symmetric group on three letters, and so $|P| = 6$.
\end{remark}


\end{document}